\documentclass[12pt, a4paper]{article}
\usepackage[english, activeacute]{babel}
\usepackage{amssymb,amsmath,latexsym,epsfig,euscript}
\usepackage{graphicx}

\def\sqr#1#2{{\vcenter{\vbox{\hrule height.#2pt\hbox{\vrule width.#2pt
 height#1pt \kern#1pt \vrule width.#2pt} \hrule height.#2pt}}}}
\def\square{\mathchoice\sqr34\sqr34\sqr{2.1}3\sqr{1.5}3}

\def\arboluno{\smallmatrix \searrow &&\swarrow \\ &\bullet &\\&\downarrow &\endsmallmatrix }
\def\arboldos{\smallmatrix \searrow &&\swarrow &&\\&\bullet &&&\\ &&\searrow &&\swarrow \\ &&&\bullet &\\
&&&\downarrow &\endsmallmatrix }
\def\arboltres{\smallmatrix &&\searrow &&\swarrow \\&&&\bullet &\\ \searrow &&\swarrow &&\\ &\bullet &&&\\
&\downarrow &&&\endsmallmatrix }
\def\arbolcuatro{\smallmatrix \searrow & \downarrow &\swarrow \\ &\bullet &\\&\downarrow &\endsmallmatrix }

\newtheorem{defn}{Definition}%[section]
\newtheorem{lem}[defn]{Lemma}
\newtheorem{prop}[defn]{Proposition}
\newtheorem{thm}[defn]{Theorem}
\newtheorem{cor}[defn]{Corollary}
\newtheorem{notn}[defn]{Notation}
\newtheorem{rem}[defn]{Remark}
\newtheorem{exmpl}[defn]{Example}

           {\vspace{3.3mm}
           \noindent{\bf #1}\it}%
           {\vspace{3.3mm}}

\def\ms {\medskip}
\def\bs {\bigskip}
\def\noi{\noindent}

\begin{document}

\centerline {\Large \bf Weak Bruhat order on the set of faces of
the permutahedra.}

\vskip 1cm

\centerline {\bf Patricia PALACIOS\footnote{Depto. de Matem\'aticas, FCEN,
Univ. de Buenos Aires. Pab. I, Ciudad Universitaria (1428) Buenos Aires, Argentina. E-mail : ppalacio@bigua.dm.uba.ar}
and Mar\' \i a O. RONCO \footnote{Depto. de Matem\'atica, CBC, Univ.
de Buenos Aires. Pab. III Ciudad Universitaria (1428) Buenos Aires, Argentina. E-mail : mronco@mate.dm.uba.ar}}
\vskip 1cm
\bs

\bs

\noi {\bf Abstract.} {\it Given a Coxeter system $(W,S)$, we
define a partial order on the set of right coclasses ${\cal
P}^{(W,S)}:=\{  W_J\circ w:\ {\rm for\ }W_J\ {\rm a\ parabolic\
subgroup\ of}\ W \} $, which extends the weak Bruhat order of the
group $W$. When $W$ is the group of permutations of $n$ elements,
we get a partial order on the set of faces of the permutahedra,
which induces a partial order on the set of planar rooted trees
with $n+1$ leaves.  We use these orders to describe associative
algebra structures on the vector spaces spanned by the faces of
permutahedra on one side, and on the vector space spanned by the
planar rooted trees, on the other one. The associative products of
both algebras may be described as the sum of three different
operations, also given by the weak Bruhat order, which give to our
examples the structure of dendriform trialgebras.} \ms

 {\bf
Introduction.}  Let $S_n$ be the group of permutations of $n$
elements. In \cite{MR}, C. Malvenuto and C. Reutenauer defined a
Hopf algebra structure on the vector space $k[S_{\infty }]:=\oplus
_{n\geq 0}k[S_n]$, where $k$ is a field and $k[S_n]$ is the group
algebra of the symmetric group.  The associative product of
$k[S_{\infty }]$ is defined using the natural action of the
Solomon algebra $Sol _{\infty }$ (cf. \cite{G-K-L-L-R-T} and
\cite{LR}) on $k[S_{\infty }]$. This algebra is widely studied in
\cite{T1}.

>From another point of view, J.-L. Loday introduced in \cite{Lo} the notion of dendriform algebra. In the same work he described
free dendriform algebras in terms of planar binary rooted trees. As a particular case, he showed that the vector space
$k[{\cal Y}_{\infty }]$, spanned by all planar binary rooted trees, has a natural structure of free dendriform algebra
generated by one element. For a fixed field $k$, there exists a canonical functor from the category of dendriform algebras
into the category of associative algebras (over $k$). The associative product of $k[{\cal Y}_{\infty }]$ is also described
in \cite{Lo} and \cite{T2}.

In \cite{LR} and \cite{LR2}, the authors used the existence of a canonical surjection from $S_n$ into the set
${\cal Y}_n$ of planar binary trees with $n+1$ leaves, and a definition of the canonical generators of the Solomon algebra in
terms of the weak Bruhat order of $S_n$, to prove that :
\begin{enumerate}\item The associative structure of the Malvenuto-Reutenauer algebra comes from a dendriform algebra
structure on $k[S_{\infty }]$, which may be defined in terms of two elementary associative operations $/$ and $\backslash $
and the weak Bruhat order.
\item The weak Bruhat order of $S_n$ induces a partial order on ${\cal Y}_n$, in such a way that the associative structure of
the free dendriform algebra $k[{\cal Y}_{\infty }]$ can also be described in terms of this order, and are induced by the
dendriform algebra structure of $k[S_{\infty }]$.
\end{enumerate}
In more recent works (cf. \cite{LR3} and \cite{LR4}), the vector space $k[{\cal T}_{\infty }]$ spanned by the set of all
planar rooted trees is described as the free object on one element for a new type of algebras, called dendriform
trialgebras.

The purpose of this paper is to:
\begin{enumerate}
\item Extend the weak Bruhat order of a Coxeter group $W$ to the set ${\cal
P}^{(W,S)}:=\{  W_J\circ w:\ {\rm for\ }W_J\ {\rm a\ parabolic\ subgroup\ of}\ W\}$.
\item Prove that this order defines a structure of dendriform trialgebra on the vector
space $k[{\cal P}_{\infty}]$ spanned by the all the posets ${\cal P}_n$ associated to the symmetric groups.
\item Show that there exists a natural surjection from the set ${\cal P}_n$ into the set ${\cal
T}_n$ of planar rooted trees with $n+1$ leaves, such that the weak Bruhat order of  ${\cal P}_n$ induces a partial order
on ${\cal T}_n$. This order defines a structure of free dendriform trialgebra on the vector space $k[T_{\infty}]$, generated by
 all planar rooted trees, which coincides with the one defined in \cite{LR3}.
\end{enumerate}

The first point is described in Section 1.

In Section 2, we give a description of the poset ${\cal P}_n$ in terms of surjective maps from
$\{ 1,\dots ,n\} $ into $\{1,\dots ,r\} $, for positive integers $n$ and $r$. This description permit us to construct
subsets $SH(n_0,\dots ,n_r)$ of ${\cal P}_n$, which generalize the well-known notion of shuffles in $S_n$. We study the
orders $\subseteq $ and $\leq _B $ in this new context, define a structure of dendriform trialgebra on $k[{\cal
P}_{\infty }]$ and show that this structure may be obtained in terms of the generalized weak Bruhat order.

Section 3 is devoted to see that the orders on ${\cal P}_n$ induce orders on the set ${\cal T}_n$ of planar rooted trees with
$n+1$ leaves and to describe the free dendriform trialgebra structure of $k[{\cal T}_{\infty }]$ in terms of the weak Bruhat
order defined in this context.
\bs

\bs

\bs

\section{Weak Bruhat order}
\ms

In this Section we define an order on the set of right coclasses ${\cal P}^{(W,S)}:=\{  W_J\circ w:\ {\rm for\ }W_J\
{\rm a\ parabolic\ \ of}\ W\}$, which differs from the canonical inclusion order and extends the weak Bruhat order
of the group $W$. For the elementary definitions and results about  Coxeter groups, Solomon algebra and weak Bruhat
order we  refer to \cite{Bo}, \cite{Hu},  \cite{Sol} and \cite{B-B-H-T}. We shall deal only with finite Coxeter systems.

\begin{notn} \begin{enumerate}\item Given a Coxeter system $(W,S)$, we denote by $l(w)$ the length of an element $w\in W$.
The element $\omega _S$ is the unique element of maximal length in $W$.
\item Given a subset $J\subseteq S$, we denote by $W_J$ the standard parabolic subgroup of
$W$ generated by $J$. The pair $(W_J,J)$ is also a Coxeter system.
\item For a subset $J$ of $S$, $X_J$ denotes the Solomon subset of elements which have no descent at $J$.
We recall that $X_J:=\lbrace w\in W \mid  l(w\circ s)>l(w),\ \forall \ s\in J\rbrace $, where $\circ $ denotes the product of
the group $W$.
\end{enumerate}
\end{notn}

Let us recall the definition of the weak Bruhat order on $W$.

\begin{defn} Let $(W,S)$ be a Coxeter system, the weak Bruhat order on $W$ is defined by
 $x \leq _B x\rq \ {\rm if}\ x\rq =y\circ x,\ {\rm with}\ l(x\rq )=l(y)+l(x ).$
\end{defn}

The group $W$ equipped with the weak Bruhat order $\leq _B$ is a partially ordered set. The minimal element of
$(W, \leq _B)$ is the identity element $1_W$ of the group, and its maximal element is $\omega _S$.

Given a subset $J\subseteq S$, there exist unique elements
$\xi _J\in X_J$ and  $\omega _J\in W_J$ such that $\omega _S=\xi _J\circ \omega _J$. It is easy to check
that $\omega _J$ is the maximal element of $(W_J,J)$, and that $\xi _J$ is the longest element of $X_J$.

The following result is proved in \cite{LR2}.

\begin{lem}\label{order2} Let $(W,S)$ be a Coxeter system and let $J\subseteq S$.
\begin{enumerate}\item The subset $X_J$ of $W$ verifies that
$X_J=\lbrace w\in W \ :\   w\leq _B\xi _J\rbrace . $
\item The subgroup $W_J$ of $W$ is the set
of elements which are smaller or equal than $\omega _J$.
\end{enumerate}
\end{lem}

\begin{defn} Let $(W,S)$ be a Coxeter system. The Coxeter poset of $(W,S)$ is the set
$P^{(W,S)}$ of right coclasses modulo the parabolic standard subgroups, ordered by the inclusion $\subseteq $.
\end{defn}

The following result follows immediately from \cite{Bo} Ch. IV, p. 37, ex. 3.
\begin{lem} The set $P^{(W,S)}$ may be described as
$P^{(W,S)}:=\lbrace W_J\circ w\ {\rm with}\ J\subseteq S\ {\rm and}\ w\in X_J^{-1}\rbrace .$
\end{lem}

For $J=\emptyset $ one has that $X_{\emptyset }=W$. So, we may think the group $W$ as embedded in $P^{(W,S)}$, via the
map $w\mapsto W_{\emptyset }\circ w$.

The family of subsets ${\cal P}_r^{(W,S)}:=\{ W_J\circ w\quad :\quad \mid J\mid =r\}$, for $0\leq r\leq \mid S\mid$, defines
a graduation on the set ${\cal P}^{(W,S)}$.

Let us define an order, different from $\subseteq $, on the set $P^{(W,S)}$.

In \cite{B-B-H-T} the authors proved that for any subset $J\subseteq S$ and $s_0\in S\setminus J$, it holds that
$X_{J\cup \lbrace s_0\rbrace }\circ X_J^{W_{ J\cup \lbrace s_0\rbrace }}=X_J,$
where $X_J^{W_{J\cup \lbrace s_0\rbrace }}$ is the set of elements in $W_{J\cup \lbrace s_0 \rbrace }$
which have no descent at $J$.

Using this equality we get the following definition

\begin{defn} Let $J\subseteq S$. For any $s\in S\setminus J$, the result above
implies that there exists an element $\alpha _{J,s}\in X_J^{W_{J\cup \lbrace s\rbrace }}= X_J\cap W_{J\cup \lbrace
s\rbrace }$ such that $\xi _J=\xi _{J\cup \lbrace s\rbrace }\circ \alpha _{J,s}$.
\end{defn}

\begin{defn}\label{wBo}({\bf Weak Bruhat order on $P^{(W,S)}$}) Let $(W,S)$ be a
Coxeter system. The weak Bruhat order on the set $P^{(W,S)}$ is the transitive relation generated by the conditions:
\begin{enumerate} \item For $J\subseteq S$, $s\in S\setminus J$ and $w\in X_{J\cup \lbrace s\rbrace}^{-1}$,
it holds that $W_J\circ w< _BW_{J\cup \lbrace s\rbrace }\circ w.$
\item For $J\subseteq S$, $s\in S\setminus J$ and $w\in X_{J\cup \lbrace s\rbrace}^{-1}$,
it holds that $W_{J\cup \lbrace s\rbrace }\circ w<_B W_J\circ (\alpha _{J,s}^{-1}\circ w).$
\end{enumerate}
\end{defn}

To check that this gives a partial order relation it suffices to see that, given a coclass $W_J\circ w$ it is impossible to
get, applying several times the relations $1$ and $2$, a sequence of coclasses that ends at $W_J\circ w$. Note that applying
$2$  we get $W_J\circ w<_B W_K\circ z$ only for $w<_B z$. So, we may restrict our proof to a sequence of coclasses obtained by
applying only the first relation. But, applying only $1$, we get that $W_J\circ w<_B W_K\circ z$ implies $J$ is strictly
contained in $K$.

The weak Bruhat order of $W$ induces an order on the subset $\lbrace W_{\emptyset }\circ w\ {\rm with}\ w\in W\rbrace
\subseteq P^{(W,S)}$. The proof of the following Lemma is straigthforward.

\begin{lem}\label{order5} Let $w< _Bz$ be two elements of $W$. The coclass
$W_{\emptyset }\circ w$ is smaller than $W_{\emptyset }\circ z$ for the weak Bruhat order.
\end{lem}

Let us point out that the result above is false if we replace $\emptyset $ by any $J\subseteq S$. As an example, consider
the symmetric group $S_4$ with the set of generators $\lbrace s_1,s_2,s_3\rbrace $, where $s_i$ is the permutation that
exchanges $i$ and $i+1$. The coclasses $W_{\lbrace s_2,s_3\rbrace }\circ 1_{S_4}$ and $W_{\lbrace s_2,s_3\rbrace }\circ s_1$
are not comparable for the weak Bruhat order.

\begin{lem}\label{order6} Let $K\subseteq S$. If $z$ and $w$ are two elements of $W$
such that $z\in X_K^{-1}$ and $\omega _K\circ z \leq _B w$, then $W_K\circ z\leq _B W_J\circ w$, for any $J$ with
$w\in X_J^{-1}$.
\end{lem}

\noi {\it Proof.} We prove the lemma by induction on the number of elements of $K$.

Suppose $\mid K\mid =0$. We have that $z\leq _B\omega _K\circ z \leq _B w$, and Lemma \ref{order5} shows that
$W_{\emptyset }\circ z \leq _B W_{\emptyset }\circ w$. Now, by applying several times the relation $1$ of Definition
\ref{wBo}, we get that $W_{\emptyset }\circ w\leq _BW_J\circ w$, and the result holds.

For $\mid K\mid \geq 1$, we have that $K=K\rq \cup \lbrace s\rbrace $ with $\mid K\rq \mid = \mid K\mid -1$. The
relation $2$ of Definition \ref{wBo} states that $W_K\circ z<_BW_{K\rq }\circ (\alpha _{K\rq ,s}^{-1}\circ z)$.

\noi From the equalities
$\omega _S=\xi _{K\rq }\circ \omega _{K\rq }=\xi _K\circ \omega _K$, and $\xi _{K\rq }=\xi _K\circ \alpha _{K\rq ,s},$
we get

\noi $\omega _K=\alpha _{K\rq ,s}\circ \omega _{K\rq }$.
Since it holds that $\omega _{K\rq }\circ \alpha _{K\rq , s}^{-1}= \omega _{K\rq }^{-1}\circ \alpha _{K\rq , s}^{-1}=
\omega _K^{-1}=\omega _K$, we get
$$\omega _{K\rq }\circ (\alpha _{K\rq ,s}^{-1}\circ z)=\omega _K\circ z < w.$$
Applying a recursive argument, $W_K\circ z<_BW_{K\rq }\circ (\alpha _{K\rq ,s}^{-1}\circ z)<_BW_J\circ w$.
\hfill $\diamondsuit $

The following Theorem describes the relationship between the orders $\subseteq $ and $\leq $ on $P^{(W,S)}$.

\begin{thm} Let $(W,S)$ be a Coxeter system. For any $J\subseteq S$ and any $w\in X_J^{-1}$, the sets
$\lbrace W_K\circ z \in P^{(W,S)}\ :\ W_K\circ z\subseteq W_J\circ w\rbrace $ and $\lbrace
W_K\circ z \in P^{(W,S)}\ :\ W_{\emptyset }\circ w\leq _BW_K\circ z\leq _BW_{\emptyset }\circ
(\omega _J \circ w)\rbrace $ are equal.
\end{thm}

\noi {\it Proof.} Observe that $W_K\circ z\subseteq W_J\circ w$ if, and only if
$K\subseteq J$ and $w\leq _Bz\leq _B\omega _J\circ w$. Since $w\leq _Bz$, we have that
$W_{\emptyset }\circ w\leq _BW_{\emptyset }\circ z\leq _B W_K\circ z$.
Now, $W_K\circ z \subseteq W_J\circ w$ implies that $\omega _K\circ z\leq _B\omega _J\circ w$. Applying Lemma
 \ref{order6}, we get that $W_K\circ z\leq _B W_{\emptyset }\circ (\omega _J\circ w)=\omega _J\circ w$.

We have proved that the set $\lbrace W_K\circ z \in P^{(W,S)}\ :\ W_K\circ z\subseteq W_J\circ w\rbrace $
is contained in $\lbrace W_K\circ z \in P^{(W,S)}\ : \ W_{\emptyset }\circ w\leq _BW_K\circ z\leq _B
W_{\emptyset }\circ (\omega _J \circ w)\rbrace $.

Suppose that $W_{\emptyset }\circ w\leq _BW_K\circ z\leq _BW_{\emptyset }\circ (\omega _J \circ w)$.
>From the first inequality it is immediate that $w\leq _Bz$, from the second one we get that $z\leq _B\omega _J\circ w$.
There exists $y_0\in W_J$ such that $z=y_0\circ w$.
To check that $W_K\circ z\subseteq W_J\circ w$ we only need to prove that $K\subseteq J$. But for any
$x\in W_K$, it holds that $W_{\emptyset }\circ w\leq_B W_K\circ (x\circ z)\leq _B W_{\emptyset}\circ (\omega _J\circ w)$,
which implies that $w\leq _B x\circ z\leq _B\omega _J\circ w$. So, $x\circ z = y_x\circ w$, for some $y_x\in W_J$;
consequently $x=y_x\circ y_0^{-1}$. We get that $W_K\subseteq W_J$, which implies (cf. \cite{Bo}) that $K\subseteq
J$. \hfill $\diamondsuit $
\bs

\bs

\bs

\section{Symmetric groups}
\ms

Let $(S_n,\circ )$ be the group of permutations of $n$ elements. The group $S_n$ is a Coxeter group generated by $n-1$
transpositions $s_1,\dots ,s_{n-1}$, where $s_i$ is the permutation that exchanges $i$ and $i+1$.
We shall denote the set $\lbrace s_1,\dots ,s_{n-1}\rbrace $ by ${\cal S}_n$.

\noi The longest element of $S_n$ is the permutation $\omega _n:=s_1\circ s_2\circ s_1\circ \dots s_{n-2}\circ \dots \circ
s_1\circ s_n\circ s-{n-1}\dots \circ s_1$.

\subsection{ The posets associated to the symmetric group}

Let $(n_1,\dots ,n_r)$ be partition of $n$ . Given elements $\sigma _i\in S_{n_i}$, for $1\leq i\leq r$, the permutation
$\sigma _1\times \dots \times \sigma _r$ in $S_n$ is defined as:
$$\sigma _1\times \dots \times \sigma _r:=(\sigma _1(1),\dots ,\sigma _1(n_1),\sigma _2(1)+n_1,\dots ,\sigma
_r(n_r)+n_1+\dots +n_{r-1}).$$
We denote by $S_{n_1,\dots ,n_r}$ the subgroup  of $S_n$ which is the image of the embedding $S_{n_1}\times \dots \times
S_{n_r}\hookrightarrow S_n$.

Let $(n_1,\dots ,n_r)$ be a partition of $n$, and let $J$ be the set
${\cal S}_n\setminus \lbrace s_{n_1},s_{n_1+n_2},\dots ,s_{n_1+\dots +n_{r-1}}\rbrace .$
It is immediate to check that the standard parabolic subgroup $W_J$ is the subgroup $S_{n_1,\dots ,n_r}$ of $S_n$.
Moreover, the longest element of $S_{n_1,\dots ,n_r}$ is $\omega _{n_1,\dots ,n_r}:=\omega _{n_1}\times \dots \times
\omega _{n_r}$.

Given a sequence of non-negative integers $n_1,\dots ,n_r$ such that $\sum\limits _{i=1}^rn_i = n$, a $(n_1,\dots
,n_r)$-shuffle is a permutation $\sigma \in S_n$ verifying that:
$$\sigma (1)<\dots <\sigma (n_1)\ ,\ \hdots \ ,
\sigma (n_1+\dots +n_{r-1}+1)<\dots <\sigma (n).$$

The set of all $(n_1,\dots ,n_r)$-shuffles in $S_n$ is denoted by $Sh(n_1,\dots ,n_r)$.

\begin{rem} \label{perm1} \begin{enumerate} \item Let $\sigma $ be an element of $S_n$ and $s_i$ a transposition in ${\cal
S}_n$. It is easily verified that $length(s_i
\circ \sigma )=length(\sigma )+1$ if, and only if, $\sigma ^{-1}(i)<\sigma ^{-1}(i+1)$.
\item Let $(n_1,\dots ,n_r)$ be a partition of $n$. The assertion above implies that the set $Sh(n_1,\dots ,n_r)$
coincides with the set $X_J$ described  in the first section, for

\noi $J={\cal S}_n\setminus \lbrace s_{n_1},s_{n_1+n_2},\dots ,s_{n_1+\dots +n_{r-1}}\rbrace .$
\end{enumerate}
\end{rem}

The longest element of $Sh(n_1,\dots ,n_r)$ (cf. \cite{LR2}) is the permutation $\xi _{n_1,\dots ,n_r}$ defined as:
$$\xi _{n_1,\dots ,n_r}(k):=k-\sum _{i=1}^{j-1}n_i+\sum _{i=j+1}^rn_i,\ {\rm for}\ \sum _{i=1}^{j-1}n_i< k\leq \sum
_{i=1}^jn_i,\ {\rm for}\ 1\leq j\leq r.$$

\begin{notn} In order to simplify notation, we denote by ${\cal P}_n$ the set $P^{(S_n,{\cal S}_n)}$.
\end{notn}

The description of the set $P^{(W,S)}$ given in Section 1 and the results above imply that for any right coclass
$W_J\circ \tau \in {\cal P}_n$, there exist a unique partition $(n_1,\dots ,n_r)$ of $n$ and a unique element $\sigma \in
Sh(n_1,\dots ,n_r)^{-1}$ such that $W_J\circ \tau =S_{n_1,\dots ,n_r}\circ \sigma $.

For any partition $(n_1,\dots ,n_r)$ of $n$ and any $1\leq k\leq r-1$, we denote by $\alpha _{n_1,\dots ,n_r}^k$ the
permutation such that $\xi _{n_1,\dots ,n_r}=\xi _{n_1,\dots ,n_k+n_{k+1},\dots ,n_r}\circ \alpha _{n_1,\dots ,n_r}^k.$

Recall that we identify an element $\sigma $ of $S_n$ with the coclass $S_{1,\dots ,1}\circ \sigma $. So we have an embedding
of partially ordered sets $(S_n,\leq _B)\hookrightarrow ({\cal P}_n,\leq _B)$.

Let us give another description of the sets ${\cal P}_n$, which makes it easier to deal with.

\begin{lem} \label{Cox1} There exists a bijection between the set ${\cal P}_n$ and the set of all surjective maps $\gamma
:\lbrace 1,\dots ,n\rbrace \rightarrow \lbrace 1,\dots ,r\rbrace $, with $r\geq 1$.
\end{lem}

\noi {\it Proof.} Let $x= S_{n_1,\dots ,n_r}\circ \sigma $ be an element of ${\cal P}_n$.

\noi Consider the map $\Psi (x):\lbrace 1,\dots ,n\rbrace \rightarrow \lbrace 1,\dots ,r\rbrace $ defined by:
$$\Psi (x)(i):=j,\ {\rm if}\ n_1+\dots +n_{j-1}<\sigma (i)\leq n_1+\dots +n_j.$$

Observe that the map $\Psi $ is bijective, the inverse map is defined by

\noi $\Psi ^{-1}(\gamma )=S_{n_1,\dots ,n_r}\circ \sigma $, where $Im(\gamma )=\lbrace 1,\dots ,r\rbrace $, $n_i$ is the number
of elements in $\gamma ^{-1}(i)$ and
$\sigma (j):=n_1+\dots +n_{\gamma (j)-1}+k$, for $\gamma ^{-1}(\gamma (j))=\lbrace i_1<\dots <i_{n_j}\rbrace $
and $ i_k=j.\hfill \diamondsuit $

Note that the image of the element $S_{1,1,\dots ,1}\circ \sigma $ under $\Psi $ is the permutation $\sigma $.

So, ${\cal P}_n=\lbrace \gamma :\lbrace 1,\dots ,n\rbrace \rightarrow \lbrace 1,\dots ,r\rbrace \quad {\rm surjective },\ {\rm
for\ some}\ r\geq 1\rbrace $.

Fix the set ${\cal P}_0:=\{ (0)\}$. Let ${\cal P}_{\infty }$ denotes the disjoint union of all the sets ${\cal P}_n$, for
$n\geq 0$. The usual composition of maps defines a multiplication on ${\cal P}_{\infty }$:
$$\gamma \circ \delta :=\left \{ \begin{array}{cc}
\gamma \circ \delta & {\rm if}\quad \gamma \in {\cal P}_r\quad {\rm and}\quad Im(\delta )=\{ 1,\dots ,r\}\\
(0)& {\rm otherwise}.\end{array}\right .$$

>From now on, we shall denote an element $\gamma \in {\cal P}_n$ by its image $\gamma =(\gamma (1),\dots ,\gamma (n))$. If
$\gamma $ and $\delta $ are permutations, then their composite coincides with their product in the symmetric group
$S_n$.

\begin{rem} Note that any element $\gamma \in {\cal P}_n$, with $n\geq 1$, may be written as $\gamma =
\gamma \rq \circ \sigma $ with $\sigma \in S_n$ and $\gamma \rq $ a non-decreasing map. The element
$\gamma \rq $ is unique, while there exist many permutations $\sigma $ verifying the equality. However there
exists a unique $\sigma $ of minimal length, which is precisely the permutation defined in Lemma \ref{Cox1}
verifying $\Psi ^{-1}(\gamma )=S_{n_1,\dots, n_r}\circ \sigma $.
\end{rem}

For $1\leq r\leq n$, the subset ${\cal P}_{n,r}:=\{ \gamma \in {\cal
P}_n\ :\  Im(\gamma )=\{ 1,\dots ,r\} \}$ corresponds to the subset ${\cal P}_r^{(S_n,{\cal S}_n}$ defined in Section 1.

There exists an embedding ${\cal P}_n\times {\cal P}_m\hookrightarrow {\cal P}_{n+m}$, given by:
$$(\gamma ,\delta )\mapsto (\gamma (1),\dots ,\gamma (n),\delta (1)+r,\dots ,\delta (m)+r),$$
where $\gamma \in {\cal P}_{n,r}$. The image of $(\gamma ,\delta )$ under this embedding is denoted by
$\gamma \times \delta $. Note that $\gamma \times (0)=\gamma =(0)\times \gamma $, for all $\gamma \in {\cal P}_{\infty }$.

For $n\geq 1$, let us denote by $t_i$ the element of ${\cal P}_{n,n-1}$ defined by:
$$t_i(j):=\left \{ \begin{array}{cc}
j & {\rm for}\quad 1\leq j\leq i\\
j-1& {\rm for}\quad i<j\leq n.\end{array}\right .$$

\begin{rem}  \label{rel} The product of the element $t_j$ with the transpositions $s_i$, which span the symmetric
group, verify the  following relationships:
$$s_i\circ t_j =\left \{ \begin{array}{cc}
t_j\circ s_i & {\rm if}\quad i<j-1\\
t_{j-1}\circ s_j\circ s_{j-1} & {\rm if}\quad i=j-1\\
t_{j+1}\circ s_j\circ s_{j+1} & {\rm if}\quad i=j\\
t_j\circ s_{i+1}& {\rm if}\quad i>j.\end{array}\right .$$
\end{rem}
\bs

\subsection{ Shuffles and wedges on ${\cal P}_n$}

We proceed to define a notion of shuffle in ${\cal P}_n$ which extends the definition of shuffle in $S_n$.

\begin{defn} Let $\gamma \in {\cal P}_n$ and $(n_1,\dots ,n_r)$ be a partition of $n$, we say that $\gamma $ is a
$(n_1,\dots ,n_r)$-shuffle in ${\cal P}_n$ if
$$\gamma (1)<\dots <\gamma (n_1),\qquad \hdots \qquad ,\gamma (n_1+\dots
+n_{r-1}+1)<\dots <\gamma (n).$$
\end{defn}

We denote by $SH(n_1,\dots ,n_r)$ the subset of $(n_1,\dots ,n_r)$-shuffles in ${\cal P}_n$. The following result is
a generalized version of  \cite{Bo}, Ch. IV, p. 37, ex. 3, to the set ${\cal P}_n$.

\begin{lem} For any element $\gamma \in {\cal P}_n$, and any $0\leq p\leq n$, there exist unique elements
$\displaystyle \gamma _1\in \bigcup_{1\leq i\leq p\atop 1\leq j\leq n-p}SH(i,j)$, $\delta _1\in {\cal P}_{p,i}$ and $\delta
_2\in {\cal P}_{n-p,j}$ such that
$\gamma = \gamma _1\circ (\delta _1\times \delta _2).$
\end{lem}

\noi {\it Proof.} Let $\gamma $ be an element of ${\cal P}_{n,r}$, there exist
unique integers  $1\leq i\leq p$ and $1\leq j\leq n-p$ and unique bijective order-preserving morphisms

\noi $\phi _1:\gamma (\{ 1,\dots ,p\})\rightarrow \{ 1,\dots ,i\}$ and $\phi _2:\gamma (\{ p+1,\dots ,n\}) \rightarrow \{
1,\dots ,j\}$.

\noi Define $\delta _1:(\phi _1(\gamma (1)),\dots ,\phi _1(\gamma (p)))$, $\delta _2:=(\phi _2(\gamma (p+1)),\dots
,\phi _2(\gamma (n)))$, and

\noi $\gamma _1:=(\phi _1^{-1}(1), \dots ,\phi _1^{-1}(i),\phi _2^{-1}(1),\dots ,\phi _2^{-1}(j))$.
It is easy to check that $\gamma _1\in SH(i,j)$ and that the elements are unique. \hfill $\diamondsuit $

\begin{defn} Let $\gamma _0,\dots ,\gamma _k$ be a family of elements of ${\cal P}_{\infty }$, with $\gamma _i\in {\cal
P}_{n_i,r_i}$ or $\gamma _i =(0)$, and let $\omega $ be an element of $SH(r_0,\dots ,r_k)$.  The wedge of
$\gamma _0,\dots ,\gamma _k$ over $\omega $ is the element of ${\cal P}_{n+k}$ given by:
$$\bigvee _{\omega }(\gamma _0,\dots ,\gamma _k):=(\omega \circ (\gamma _0\times \dots \times \gamma _k)\times 1_1)\circ
z(n_0,\dots ,n_k),$$
where $z(n_0,\dots ,n_k):=(1,\dots ,n_0,n+1,n_0+1,\dots ,n_0+n_1,n+1,\dots ,n+1, n_0+\dots +n_{k-1}+1,\dots ,n)$ and $n:=\sum
_{i=0}^kn_i$.
\end{defn}

For instance, observe that $1^n:=(1,1,\dots ,1):\{ 1,\dots ,m\} \rightarrow \{ 1\}$ is the wedge $\bigvee _{(0)}((0),\dots
,(0))$.

\begin{prop}\label{wedge} Given $\gamma \in {\cal P}_{n,r}$ there exist unique elements $\gamma _i$ in
${\cal P}_{n_i,r_i}$, with

\noi $k+\sum _{i=0}^kn_i=n$ and $r-1=\sum _{i=0}^kr_i$, and $\omega \in SH(r_0,\dots ,r_k)$
such that

\noi $\gamma =\bigvee _{\omega }(\gamma _0,\dots ,\gamma _k)$.
\end{prop}

\noi {\it Proof.}\ Let $\gamma ^{-1}(r)=\{ j_1<\dots <j_k\} $. Define the integers $n_i$, for $0\leq i\leq k$, as follows:
$$n_0:=j_1-1,\qquad \dots ,\qquad n_l:=j_{l+1}-j_l-1,\qquad \dots ,\qquad n_k:=n-j_k.$$
It is easy to check that $\gamma = (\alpha \times 1_1)\circ z(n_0,\dots ,n_k)$, with $\alpha \in {\cal P}_{n-k,r-1}$.
>From Lemma \ref{Cox1}, we get that there exist unique $\gamma _0\in {\cal P}_{n_0},\dots ,\gamma _k\in {\cal P}_{n_k}$
 and $\omega \in SH(r_0,\dots ,r_k)$ such that $\alpha =\omega \circ (\gamma _0\times \dots \gamma _k)$.

\noi The argument above proves the existence of the decomposition. It is clear that $\alpha $ and $n_0,\dots ,n_k$ are
unique. Lemma \ref{Cox1} implies the unicity of the $\gamma _i$\rq and $\omega $. \hfill $\diamondsuit $

\begin{defn} Let $\omega \in SH(r_0,\dots ,r_k)$ and let  $1\leq j\leq \sum _{i=0}^kr_i$. We associate to $\omega $ and $j$
a family of positive integers $(s_0^j,\dots ,s_k^j)$, with $s_i^j=r_i$ or $s_i^j=r_i-1$, and a shuffle
$\omega _j\in SH(s_0^j,\dots ,s_k^j)$ in the following way:
\begin{enumerate}\item If $t_j\circ w\in SH(r_0,\dots ,r_k)$, then $(s_0^j,\dots ,s_k^j):=(r_0,\dots ,r_k)$ and $\omega
_j:=t_j\circ \omega $.
\item  If $t_j\circ \omega \notin SH(r_0,\dots ,r_k)$, then
$$s_i^j:=\left \{ \begin{array}{cc}
r_i & {\rm if}\quad \{ j,j+1\}\not\subseteq \omega (\{ R_{i-1}+1, \dots ,R_i\} ), \\
r_i-1& {\rm if}\quad \{ j,j+1\}\subseteq \omega (\{ R_{i-1}+1, \dots ,R_i\} ), \end{array}\right .$$
where $R_l:=r_0+\dots +r_l$.
The element $\omega _j$ is the $(s_0^j,\dots ,s_k^j)$-shuffle in ${\cal P}_{\sum s_i^j}$ such that, for $0\leq i\leq k$ and
$1\leq l\leq s_i^j$, $\omega _j(s_0^j+\dots +s_{i-1}^j+l):=$
$$\left \{ \begin{array}{cc}
\omega (R_{i-1}+l)& {\rm if}\quad s_i^j=r_i,\\
\omega (R_{i-1}+l)& {\rm if}\quad s_i^j=r_i-1\quad {\rm and}\quad \omega (R_{i-1}+l)\leq j,\\
\omega (R_{i-1}+l+1)-1&{\rm if}\quad s_i^j=r_i-1\quad {\rm and}\quad \omega (R_{i-1}+l)>j.\end{array}\right .$$
\end{enumerate}
\end{defn}

\begin{lem}\label{prodwedge} Let $\omega \in SH(r_0,\dots ,r_k)$ and $\gamma _i\in {\cal P}_{n_i,r_i}$, for $0\leq i\leq k$.

\noi For $1\leq j\leq r:=\sum _{i=0}^kr_i+1$, it holds that:
\begin{enumerate}\item If $j< r-1$ and $t_j\circ \omega \in SH(r_0,\dots ,r_k)$, then
$t_j\circ \bigvee _{\omega }(\gamma _0,\dots ,\gamma _k)=\bigvee _{t_j\circ \omega }(\gamma _0,\dots ,\gamma _k).$
\item If $j< r-1$ and $t_j\circ \omega \notin SH(r_0,\dots ,r_k)$, then
$t_j\circ \bigvee _{\omega }(\gamma _0,\dots ,\gamma _k)=\bigvee _{\omega _j}(\gamma _0^j,\dots ,\gamma _k^j),$
where $\gamma _i^j:=\left \{ \begin{array}{cc}
\gamma _i&\quad {\rm if}\quad s_i^j=r_i,\\
t_{j_i} \circ \gamma _i&\quad {\rm if}\quad s_i^j=r_i-1,\end{array}\right .$

\noi for $\omega ^{-1}(j)=r_0+\dots +r_{i-1}+j_i$, with $0\leq i\leq k$.
\end{enumerate}
\end{lem}

\noi {\it Proof.} The proof follows straightforward, using the definition above. \hfill $\diamondsuit $
\bs

\subsection{ Inclusion order and Weak Bruhat order on ${\cal P}_n$}

Translating the definition of inclusion order and weak Bruhat order to the set of surjective maps $\gamma :\{ 1,\dots
,n\}\rightarrow \{ 1,\dots ,r\}$, we get the following result.

\begin{lem} \begin{enumerate}\item Let $\gamma \in {\cal P}_{n,r}$ and $\delta \in {\cal P}_{n,s}$. It holds that
$\gamma  \subseteq \delta $ in ${\cal P}_n$ if, and only if, there exists a non-decreasing map $\rho \in {\cal P}_{r,s}$ such
that $\delta =\rho \circ \gamma $.
\item The weak Bruhat order $\leq _B$ on the set ${\cal P}_n$ is the transitive relation spanned by the conditions:
\begin{enumerate}
\item If $\gamma ^{-1}(i)<\gamma ^{-1}(i+1)$ for some $1\leq i\leq r$ then  $\gamma <_B
t_i\circ \gamma ,$
\item If $\gamma ^{-1}(i) > \gamma ^{-1}(i+1)$ for some $1\leq i\leq r$ then $t_i\circ \gamma
<_B \gamma ,$
\end{enumerate}
for $\gamma \in {\cal P}_{n,r} $, where $\gamma ^{-1}(i)<\gamma ^{-1}(i+1)$ means that any element $j$ of $\gamma ^{-1}(i)$
verifies $j<k$ for all $k\in \gamma ^{-1}(i+1)$.
\end{enumerate}
\end{lem}

\noi{\it Proof.} The proof of the first point is immediate.

For the second one, it is easily seen that the weak Bruhat order $\leq _B$ defined on the elements of ${\cal P}_n$, seen as
right coclasses, is the transitive relation spanned by the following conditions:
\begin{enumerate}
\item $S_{n_1,\dots ,n_r}\circ \sigma <_BS_{n_1,\dots ,n_k+n_{k+1},\dots ,n_r}\circ \sigma .$
\item $S_{n_1,\dots ,n_k+n_{k+1},\dots ,n_r}\circ \sigma <_BS_{n_1,\dots ,n_r}\circ (\alpha _{n_1,\dots ,n_r}^k)^{-1}\circ
\sigma .$
\end{enumerate}
for any partition $(n_1,\dots ,n_r)$ of $n$, any $1\leq k\leq r-1$ and any $\sigma \in Sh(n_1,\dots ,n_k+n_{k+1},\dots
,n_r)$.
Translating these conditions from the elements of type  $S_{n_1,\dots ,n_r}\circ \sigma $ to the surjective maps $\{ 1,\dots
,n\} \rightarrow \{1,\dots ,r\} $ the assertion follows. \hfill $\diamondsuit $

\begin{rem} \begin{enumerate} \item Given $\gamma \in {\cal P}_n$ it holds that $1_n\leq _B\gamma \leq _B \omega _n$, where
$1_n$ is the identity element of $S_n$ and $\omega _n = (n,n-1,\dots ,1)$.
\item Suppose $\gamma =\tau \circ \sigma $, with $\tau $ a non-decreasing map and $\sigma \in S_n$. If

\noi $\gamma \in SH(n_1,\dots ,n_r)$, then $\sigma \in Sh(n_1,\dots ,n_r)$.
\item Suppose $\gamma \leq _B \gamma \rq $ in ${\cal P}_n$. If $\gamma \rq (j)<\gamma \rq (k)$ for some $1\leq j<k\leq n$,
then $\gamma (j)<\gamma (k)$. Conversely, if $\gamma (j)>\gamma (k)$ for some $1\leq j<k\leq n$,
then $\gamma \rq (j)>\gamma \rq (k)$.
\end{enumerate}
\end{rem}

\begin{prop} \label{shuf} It holds that
$$SH(n_1,\dots ,n_r)=\{ \omega \in {\cal P}_n\quad {\rm such\quad
that}\quad \omega \leq _B \xi _{n_1,\dots ,n_r}\}.$$
\end{prop}

\noi {\it Proof.} Suppose $\omega \in {\cal P}_n$ is such that $\omega \leq _B \xi _{n_1,\dots ,n_r}$. The Remark above
implies that $\omega \in SH(n_1,\dots ,n_r)$.

Conversely, suppose $\omega \in SH(n_1,\dots ,n_r)$. We have that $\omega =\tau \circ \sigma $, with $\tau $ a non-decreasing
map and $\sigma \in Sh(n_1,\dots ,n_r)$, which may be supposed of minimal length. But $\tau $ may be written as a product
$\tau = t_{j_k}\circ \dots \circ  t_{j_1}$, with $j_1\leq j_2\leq \dots \leq j_k$, for unique integers $j_1,\dots
,j_k$. If $k=0$, the result follows from the Remark above.

\noi For $k=1$, we have two possibilities:
\begin{enumerate}
\item If $\sigma ^{-1}(j_1+1)<\sigma ^{-1}(j_1)$, then $\omega <_B\sigma \leq _B \xi _{n_1,\dots ,n_r}$.
\item If $\sigma ^{-1}(j_1)<\sigma ^{-1}(j_1+1)$, then $\omega <_B s_{j_1}\circ \sigma $. Since $\omega \in SH(n_1,\dots
,n_r)$, it holds that there exists $0\leq l\leq r-1$ such that $\sigma ^{-1}(j_1)\leq n_1+\dots +n_l <\sigma ^{-1}(j_1+1)$,
which implies that $s_{j_1}\circ \sigma \in Sh(n_1,\dots ,n_r)$. So, $\omega <s_{j_1}\circ \sigma \leq \xi _{n_1,\dots ,n_r}$.
\end{enumerate}

\noi For $k>1$, let $\gamma :=t_{j_{k-1}}\circ \hdots \circ  t_{j_1}\circ \sigma$. The element $s_{j_k}\circ \gamma $
verifies that $\omega <_B s_{j_k}\circ \gamma $, $s_{j_k}\circ \gamma \in SH(n_1,\dots ,n_r)$ and, using Remark \ref{rel}, we
get that
$$s_{j_k}\circ \gamma =t_{l_{k-1}}\circ \hdots \circ t_{l_1}\circ \sigma \rq,\quad {\rm for\quad some}\quad \sigma \rq
\in S_n.$$
So, applying recursive hypothesis to $s_{j_k}\circ \gamma $, we get that $\omega <_Bs_{j_k}\circ \gamma \leq _B \xi _{n_1,\dots
,n_r}$. \hfill $\diamondsuit $

\begin{prop} Let $\gamma _0,\dots ,\gamma _k$ and $\delta _0,\dots ,\delta _k$ be families of elements in ${\cal P}_{\infty }$
verifying that $\gamma _i\leq _B \delta _i$ in ${\cal P}_{n_i}$. It holds that
$\gamma _0\times \dots \times \gamma _k\leq _B\delta _0\times \dots \times \delta _k.$
\end{prop}

\noi {\it Proof.} Clearly, it suffices to show that
$$\gamma _0\times \dots \times \gamma _k\leq _B \gamma _0\times \dots \gamma _{i-1}\times \delta _i\times \gamma
_{i+1}\times \dots \times \gamma _k,$$
for all $0\leq i\leq k$.

Fix $0\leq i\leq k$, and suppose $\gamma _l\in {\cal P}_{n_l,r_l}$, for $0\leq l\leq k$. The proof may be restricted to the
following two cases:
\begin{enumerate} \item When $\gamma _i^{-1}(j)<\gamma _i^{-1}(j+1)$, for some $1\leq j\leq r_i$, and $\delta _i=t_j\circ
\gamma _i$.
\item When $\delta _i^{-1}(j)>\delta _i^{-1}(j+1)$, for some $j$, and $\gamma _i=t_j\circ \delta _i$.
\end{enumerate}
The proofs of both cases are analogous, so we give the proof of the first one.

If $\gamma _i^{-1}(j)<\gamma _i^{-1}(j+1)$, for some $1\leq j\leq r_i$ and $\delta _i=t_j\circ \gamma
_i$, then
$$(\gamma _0\times \dots \gamma _k)^{-1}(r_0+\dots +r_{i-1}+j)<(\gamma _0\times \dots \gamma _k)^{-1}(R_{i-1}+j+1),\quad {\rm
and}$$
$$\gamma _0\times \dots \times \gamma _{i-1}\times \delta _i\times \gamma _{i+1}\times \dots \times \gamma _k=
t_{R_{i-1}+j}\circ (\gamma _0\times \dots \times \gamma _k),$$ where $R_{i-1}:=r_0+\dots +r_{i-1}$.

\noi So, $\gamma _0\times \dots \times \gamma _k<_B \gamma _0\times \dots \times \gamma _{i-1}\times \delta _i\times
\gamma _{i+1}\times \dots \times \gamma _k$. \hfill $\diamondsuit $

\begin{lem}\label{shufs} Let $\gamma _0,\dots ,\gamma _k$ be a family of elements in ${\cal P}_{\infty }$ such that $\gamma
_i\in {\cal P}_{n_i,r_i}$, and let $\omega $ and $\omega \rq $ be elements of $SH(r_0,\dots ,r_k)$ such that $\omega \leq
_B\omega
\rq$. It holds that:
$$\omega \circ (\gamma _0\times \dots \times \gamma _k)\leq _B\omega \rq \circ (\gamma _0\times \dots \times \gamma _k).$$
\end{lem}

\noi {\it Proof.} Again, it suffices to consider the cases:
\begin{enumerate}
\item $\omega ^{-1}(j)<\omega ^{-1}(j+1)$ for some $j$, and $\omega \rq =t_j\circ \omega $.
\item $(\omega \rq )^{-1}(j)>(\omega \rq )^{-1}(j+1)$ for some $j$, and $\omega =t_j\circ \omega \rq$.
\end{enumerate}

Since the proofs of the cases are very similar, we give the proof of the second one.

Suppose $(\omega \rq )^{-1}(j)>(\omega \rq )^{-1}(j+1)$ for some $j$, and $\omega =t_j\circ \omega \rq$. Since $\omega
\in SH(r_0,\dots ,r_k)$, there must exists $0\leq l\leq k-1$ such that
$(\omega \rq )^{-1}(j+1)\subseteq \{ 1,\dots ,n_0+\dots +n_l\} $ and $(\omega \rq )^{-1}(j)\subseteq \{n_0+\dots +n_l+1,\dots
,n_0+\dots +n_k\} $. So, we get that
$$(\omega \rq \circ (\gamma _0\times \dots \times \gamma _k))^{-1}(j)>(\omega \rq \circ (\gamma _0\times \dots \times \gamma
_k))^{-1}(j+1),$$
and $\omega \circ (\gamma _0\times \dots \times \gamma _k)=t_j\circ \omega \rq \circ (\gamma _0\times \dots \times \gamma
_k).$
The argument above implies that:
$$\omega \circ  (\gamma _0\times \dots \times \gamma _k)\leq _B\omega \rq  (\gamma _0\times \dots \times
\gamma _k). \qquad \qquad \hfill \diamondsuit $$

\begin{cor}\label{wBor1} Let $\gamma _0,\dots ,\gamma _k$ be elements in ${\cal P}_{\infty }$ with $\gamma _i\in {\cal
P}_{n_i,r_i}$. If $\omega \leq _B\omega \rq $ in $SH(r_0,\dots ,r_k)$, then
$\bigvee _{\omega }(\gamma _0,\dots ,\gamma _k)\leq _B\bigvee _{\omega \rq }(\gamma _0,\dots ,\gamma _k).$
\end{cor}

\begin{lem}\label{wBor2} Let $\gamma _i\in {\cal P}_{n_i,r_i}$, for $0\leq i\leq k$, and let $\omega
\in SH(r_0,\dots ,r_k)\cap S_{r-1}=Sh(r_0,\dots ,r_k)$.
\begin{enumerate} \item If $\gamma _l\leq _B t_j\circ \gamma _l$ and $\omega (r_0+\dots +r_{l-1}+j+1)=\omega (r_0+\dots +
r_{l-1}+j)+1$, for some $0\leq l\leq k$ and $1\leq j\leq r_l$, then:
$$\bigvee _{\omega }(\gamma _0,\dots ,\gamma _k)\leq _B\bigvee _{\omega \rq }(\gamma _0,\dots ,t_j\circ \gamma
_l,\dots ,\gamma _k),$$
where $\omega \rq :=t_{\omega (r_0+\dots +r_{l-1}+j)}\circ \omega $.
\item If $t_j\circ \gamma _l\leq _B  \gamma _l$ and $\omega (r_0+\dots +r_{l-1}+j+1)=\omega (r_0+\dots +
r_{l-1}+j)+1$, for some $0\leq l\leq k$ and $1\leq j\leq r_l$, then:
$$\bigvee _{\omega \rq}(\gamma _0,\dots ,t_j\circ \gamma _l,\dots ,\gamma _k)\leq _B\bigvee _{\omega }(\gamma _0,\dots
,\gamma _k),$$
where $\omega \rq :=t_{\omega (r_0+\dots +r_{l-1}+j)}\circ \omega $.
\end{enumerate}
\end{lem}

\noi {\it Proof.} The proof is immediate, using that $t_{\omega (r_0+\dots +r_{l-1}+j)}\circ \omega \in SH(r_0,\dots
,r_k)$.\hfill $\diamondsuit $
\bs

\subsection {Dendriform trialgebra structure on the vector space spanned by ${\cal P}_{\infty }$}

We want to extend some well-known results for the shuffles in $S_n$ to our generalized shuffles in ${\cal P}_n$. Let us begin
with the associativity of the shuffle.

\begin{prop} ({\bf Associativity of $SH$}) Given non-negative integers $n,m,r$ the following equality holds:
$$\displaylines {
SH(n,m,r)=\bigcup _{1\leq j\leq n+m}SH(j,r)\circ (SH(n,m)\times 1_r)=\bigcup _{1\leq k\leq m+r}SH(n,k)\circ (1_n\times
SH(m,r))\cr }$$
\end{prop}

\noi {\it Proof.} From Lemma \ref{Cox1}, we get that for any $\gamma \in SH(n,m,r)$ there exist unique integers $1\leq j\leq
n+m$ and $1\leq k\leq m+r$, and unique elements $\gamma _1\in SH(j,r)$, $\gamma _2\in SH(n,k)$, $\delta _1\in {\cal P}_{n+m}$,
$\delta _2\in {\cal P}_r$, $\tau _1\in {\cal P}_n$ and $\tau _2\in {\cal P}_{m+r}$, such that
$$\gamma = \gamma _1\circ (\delta _1\times \delta _2)=\gamma _2\circ (\tau _1\times \tau _2).$$
To end the proof note that:

\noi $(i)$ Since $\gamma \in SH(n,m,r)$ and $\gamma _1\in SH(j,r)$,the element $\delta _1$ must belong to
$SH(n,m)$ and $\delta _2$ must be equal to $1_r$.

\noi $(ii)$ The fact that $\gamma \in SH(n,m,r)$ and that $\gamma _2\in SH(n,k)$ imply that $\tau _1=1_n$ and $\tau _2\in
SH(m,r)$. \hfill $\diamondsuit $

For $n,m\geq 0$, the set $SH(n,m)$ is the disjoint union of the following three subsets:
$$\displaylines {
SH^{\succ }(n,m):=\{ \gamma \in SH(n,m):\quad \gamma (n)<\gamma (n+m)\},\cr
SH^{\bullet }(n,m):=\{ \gamma \in SH(n,m):\quad \gamma (n)=\gamma (n+m)\},\cr
SH^{\prec }(n,m):=\{ \gamma \in SH(n,m):\quad \gamma (n)>\gamma (n+m)\}.\cr }$$

\begin{lem} \label{wbo} For $n,m\geq 0$, the sets $SH^{\succ }(n,m)$, $SH^{\bullet }(n,m)$ and $SH^{\prec }(n,m)$ may be
described in terms of the weak Bruhat order as follows:
\begin{enumerate} \item $SH^{\succ }(n,m)=\{ \gamma \in {\cal P}_{n+m}:\gamma \leq _B \xi
_{n,m-1}\times 1_1\} ,$
\item $SH^{\bullet }(n,m)=\{  \gamma \in {\cal P}_{n+m}:z(n-1,m-1,0)\leq _B\gamma \leq _B
(\xi _{n-1,m-1}\times 1_1)\circ z(n-1,m-1,0)\} ,$
\item $SH^{\prec }(n,m)=\{   \gamma \in {\cal P}_{n+m}:z(n-1,m)\leq _B \gamma \leq _B \xi
_{n,m}\} ,$
\end{enumerate}
where $z(n-1,m-1,0)=(1,\dots ,n-1,n+m-1,n,\dots ,n+m-1)$ and $z(n-1,m)=(1,\dots ,n-1,n+m,n,\dots ,n+m-1)$.
\end{lem}

\noi {\it Proof.} The proof of the formulas follows from Proposition \ref{shuf}. We give the proof of the first
equality, the others may be checked in a similar way.

If $\gamma \in SH^{\succ }(n,m)$, then $\gamma = \gamma _1\times 1_1$, with $\gamma _1\in SH(n,m-1)$. By Propositition
\ref{shuf}, we get that $\gamma _1\leq _B \xi _{n,m-1}$. It implies that $\gamma \leq _B \xi _{n,m-1}\times 1_1$.

Conversely, let $\gamma \in {\cal P}_{n+m}$ be such that $\gamma \leq _B \xi _{n,m-1}\times 1_1$. It is easily seen that
$\gamma =\gamma _1\times 1_1$, with $\gamma _1\leq _B \xi _{n,m-1}$. So, $\gamma \in SH^{\succ }(n,m)$. \hfill $\diamondsuit $

\begin{rem}\label{trish} It is not difficult to check that the following identities hold:
\begin{enumerate}
\item $\bigcup _{i=1}^{m+r}SH^{\succ }(n,i)\circ (1_n\times SH^{\succ }(m,r))=\bigcup _{j=1}^{n+m}SH^{\succ }(j,r)\circ
(SH(n,m)\times 1_r),$
\item $\bigcup _{i=1}^{m+r}SH^{\succ }(n.i)\circ (1_n\times SH^{\prec }(m,r))=\bigcup _{j=1}^{n+m}SH^{\prec }(j,r)\circ
(SH^{\succ }(n,m)\times 1_r),$
\item $\bigcup _{i=1}^{m+r}SH^{\prec }(n,i)\circ (1_n\times SH(m,r))=\bigcup _{j=1}^{n+m}SH^{\prec }(j,r)\circ (SH^{\prec
}(n,m)\times 1_r),$
\item $\bigcup _{i=1}^{m+r}SH^{\bullet }(n,i)\circ (1_n\times SH^{\cdot }(m,r))=\bigcup _{j=1}^{n+m}SH^{\bullet }(j,r)\circ
(SH^{\bullet }(n,m)\times 1_r),$
\item $\bigcup _{i=1}^{m+r}SH^{\succ }(n,i)\circ (1_n\times SH^{\bullet }(m,r))=\bigcup _{j=1}^{n+m}SH^{\bullet }(j,r)\circ
(SH^{\succ }(n,m)\times 1_r),$
\item $\bigcup _{i=1}^{m+r}SH^{\bullet }(n,i)\circ (1_n\times SH^{\succ }(m,r))=\bigcup _{j=1}^{n+m}SH^{\bullet }(j,r)\circ
(SH^{\prec }(n,m)\times 1_r),$
\item $\bigcup _{i=1}^{m+r}SH^{\bullet }(n,i)\circ (1_n\times SH^{\prec}(m.r))=\bigcup _{j=1}^{n+m}SH^{\prec }(j,r)\circ
(SH^{\bullet }(n,m)\times 1_r).$
\end{enumerate}
\end{rem}
\bs

The following definition is given in \cite{LR3}, it generalizes the definition of dendriform algebra introduced by J.-L. Loday
in \cite{Lo}.

\begin{defn} Let $k$ be a field. A dendriform trialgebra over $k$ is a vector space $V$ equipped with three $k$-linear maps
$\succ \ ,\  \cdot \ ,\ \prec \ :V\otimes _kV\rightarrow V$ verifying the following identities:
\begin{enumerate}\item $v\succ (w\succ u)=(v*w)\succ u,$
\item $v\succ (w\prec u)=(v\succ w)\prec u,$
\item $v\prec (w*u)=(v\prec w)\prec u,$
\item $v\cdot (w\cdot u)=(v\cdot w)\cdot u,$
\item $v\succ (w\cdot u)=(v\succ w)\cdot u,$
\item $v\cdot (w\succ u)=(v\prec w)\cdot u,$
\item $v\cdot (w\prec u)=(v\cdot w)\prec u,$
\end{enumerate}
for all $v,w$ and $u$ in $V$, where $v*w:=v\succ w+v\cdot w+v\prec w$.
\end{defn}

If $(V,\succ ,\cdot ,\prec )$ is a tridendriform algebra, then $(V,*)$ is an associative algebra over $k$.

Let $k[{\cal P}_{\inftyÊ}]$ be the free $k$ vector space spanned by the set ${\cal P}_{\infty }$. We define operations $\succ
$, $\cdot $ and $\prec $ on $k[{\cal P}_{\infty }]$ as follows:
\begin{enumerate}\item $\gamma \succ \delta :=\sum _{\omega \in SH^{\succ }(r,s)}\omega \circ (\gamma \times \delta )$,
\item $\gamma \cdot \delta :=\sum _{\omega \in SH^{\bullet }(r,s)}\omega \circ (\gamma \times \delta )$,
\item $\gamma \prec \delta :=\sum _{\omega \in SH^{\prec }(r,s)}\omega \circ (\gamma \times \delta )$,
\end{enumerate}
for $\gamma \in {\cal P}_{n,r}$ and $\delta \in {\cal P}_{m,s}$.

Remark \ref{trish} implies that $k[{\cal P}_{\infty }]$ is a dendriform trialgebra over $k$. Moreover, it is easy to see that
the associative product $*$ is given by the formula:
$$\gamma *\delta = \sum _{\omega \in SH(r,s)}\omega \circ (\gamma \times \delta ),\quad {\rm for}\quad \gamma \in {\cal
P}_{n,r}\quad {\rm and}\quad \delta \in {\cal P}_{m,s}.$$

\begin{thm} The operations of the dendriform trialgebra $k[{\cal P}_{\infty }]$ may be described in terms of the weak
Bruhat order as follows:
$$\displaylines {
\gamma \succ \delta =\sum _{\gamma \times \delta \leq _B \omega \leq _B (\xi _{r,s-1}\times 1_1)\circ (\gamma \times \delta )}
\omega ,\cr
\gamma \cdot \delta =\sum _{z(r-1,s-1,0)\circ (\gamma \times \delta )\leq _B\omega \leq _B (\xi _{r-1,s-1}\times
1_1)\circ z(r-1,s-1,0)\circ (\gamma \times \delta )} \omega ,\cr
\gamma \prec \delta =\sum _{z(r-1,s)\circ (\gamma \times \delta )\leq _B\omega \leq _B \xi _{r,s}(\gamma \times \delta
)}\omega ,\cr
\gamma *\delta =\sum _{\gamma \times \delta \leq _B \omega \leq _B \xi _{r,s}\circ (\gamma \times \delta }\omega ,\cr }$$
where $\gamma \in {\cal P}_{n,r}$, $\delta \in {\cal P}_{m,s}$.
\end{thm}

\noi {\it Proof.} It suffices to apply Lemmas \ref{shufs} and \ref{wbo}. \hfill $\diamondsuit $
\bs

\bs

\bs

\section{Planar rooted trees}
\ms

In this paper, a {\it tree} will always be a finite planar non-empty oriented connected graph $t$ without loops, and such
that for any vertex of $t$ there are at least two incoming edges and exactly one outgoing edge.

For $n \geq 0$, let ${\cal T}_n$ denote the set of planar rooted trees with $n +1$ leaves:

${\cal T}_0=\{\downarrow \}$, ${\cal T}_1=\{ \arboluno \}$, ${\cal T}_2=\{ \arboldos ,
\arbolcuatro, \arboltres\}.$

\noi The {\it degree} of a tree $t$ is $n$ when $t$ belongs to ${\cal T}_n$. We denote the degree of $t$ by $\mid t\mid $.

A tree $t$ is called {\it binary} if each vertex of $t$ has exactly two incoming edges. We denote by ${\cal Y}_n$ the subset
of binary trees of ${\cal T}_n$. Note that the elements of ${\cal Y}_n$ are the trees which have exactly $n+1$ leaves and $n$
vertices.

The tree $\varsigma _n$ is the element of ${\cal T}_n$ which has exactly one vertex.

The {\it wedge} of planar rooted trees $t^0 \in {\cal T}_{n_0}, \ldots, t^r\in {\cal T}_{n_r}$, is the tree $\bigvee(t^0,
\ldots, t^r)$ of degree $\sum _{i=0}^rn_i + r-1$, obtained by joining the
roots of $t^0, \ldots, t^r$ to a new vertex and creating a new root. For any tree $t \not = \downarrow$, there exist unique
trees
$t^0,
\ldots, t^r$ such that $t= \bigvee(t^0, \ldots, t^r)$.

\noi We have that ${\cal T}_n=\bigsqcup\limits_{n_0 + \ldots + n_r + r -1=n} {\cal T}_{n_1} \times \ldots \times
{\cal T}_{n_r}$, for $n \geq 1$.

We denote by ${\cal T}_{\infty }$ the graded set of all planar trees $\bigcup _{n\geq 0} {\cal T}_n$.
\bs

\subsection{Orders on ${\cal T}_n$}

We are going to define two orders $\subseteq $ and $\leq _B$ on the set of planar trees ${\cal T}_n$. The first one is the
opposite of the order described by V. Ginzburg and M. Kapranov in \cite{GK}. In the next section we show that both orders are
induced by the orders $\subseteq $ and $\leq _B$ on the set ${\cal P}_n$.

Let $t$ be a tree and let $e$ be an internal edge of $t$. We denote by $t_e$ the tree obtained from $t$ by contracting
$e$ into a point. For instance, if $t=\smallmatrix \searrow &&\swarrow &&\searrow&&\swarrow\\&\bullet &&&&\bullet &\\&&\searrow
&&\swarrow e&&\\&&&\bullet &&&\\&&&\downarrow&&&\endsmallmatrix $, then $t_e=\smallmatrix \searrow &&\swarrow &&\\&\bullet &&&\\ &&\searrow
&\downarrow&\swarrow \\ &&&\bullet &\\ &&&\downarrow &\endsmallmatrix $.

\begin{defn} Let $\subseteq $ be the order on ${\cal T}_n$ transitively generated by the relation:
$$t\subseteq t_e,\ {\rm for\ any\ internal\ edge}\ e\ {\rm of}\ t.$$
\end{defn}

Note that the binary trees of ${\cal T}_n$ are the minimal elements for $\subseteq $, while the element $\varsigma _n$ is the
maximal one.

\begin{defn}\label{BOT} Let $\leq _B$ be the relation on ${\cal T}_n$ generated transitively by the following relations:
\begin{enumerate}
\item If $t^{i_0} < _B w^{i_0} \in {\cal T}_{n_{i_0}}$ , then
$\bigvee(t^0, \ldots , t^{i_0},\ldots , t^k) < _B\bigvee(t^0, \ldots , w^{i_0},\ldots , t^k)$ in
${\cal T}_{n_0 + \ldots + n_k + k -1}$,
\item If $t=\bigvee(t^0, \ldots, t^k)\in {\cal T}_n$ and $w=\bigvee(w^0,
\ldots, w^h) \in {\cal T}_m$, then
$$\bigvee(t,w^0, \ldots, w^h) < _B\bigvee(t^0, \ldots, t^k, w^0, \ldots, w^h)$$
\item If $t=\bigvee(t^0, \ldots, t^k)\in {\cal T}_n$, then
$$t < _B\bigvee(t^0, \ldots, t^j, \bigvee(t^{j + 1}, \ldots, t^k)),$$
for $1 \leq j < k$.
\end{enumerate}
\end{defn}

\begin{exmpl}
On ${\cal T}_2$, we have
$$ \arboldos \leq _B\arbolcuatro \leq _B\arboltres $$
\end{exmpl}

\begin{prop} The pair $({\cal T}_n, \leq _B)$ is a partially ordered set, for all $n\geq 1$.
\end{prop}

\noi {\it Proof.} We need to check that, given a tree $t$ in ${\cal T}_n$, it does not exist a sequence of trees
$t=t_0,t_1,\dots,t_m=t$ such that $t_{i+1}$ is obtained from $t_i$ applying one of the relations $1,\ 2$ or $3$ of
Definition \ref{BOT}.

\noi We proceed by induction on $n$. For $n=1$ or $2$ the result is immediate. Suppose now that $\mid t\mid > 2$,
there exist trees $t^0,\dots ,t^k$ such that $t=\bigvee (t^0,\dots ,t^k)$. Note that, if $t_i$ is obtained from $t_{i-1}$ by
 an operation of type $1$ or $3$, then $\mid t_i^1\mid = \mid t_{i-1}^1\mid $, while applying operation of type $2$, we get
$\mid t_i^1\mid <\mid t_{i-1}^1\mid $. So, it there exists a sequence $t=t_0,t_1,\dots,t_m=t$ such that $t_{i+1}$ it must be
obtained applying only operations of types $1$ or $3$.

\noi Now, suppose $t_j=\bigvee (t_j^0,\dots ,t_j^{k_j})$, for $0\leq j\leq m$. If $t_j$ is obtained from $t_{j-1}$ by
an operation of type $1$, then $\mid t_j^{r_j}\mid = \mid t_{j-1}^{r_{j-1}}\mid $, and if it is obtained from $t_{j-1}$ by
an operation of type $2$, then $\mid t_j^{r_j}\mid > \mid t_{j-1}^{r_{j-1}}\mid $. We may conclude then that the sequence
$t=t_0,t_1,\dots,t_m=t$ must be obtained applying only operations of type $1$.

\noi In this case, it must exist at least one $1\leq k\leq m$ and a sequence $t^k=t_0^k<_B\dots <_Bt_m^k=t^k$. But,
by induction hypothesis, we know that it is false. So, the relation $\leq _B$ is a well defined order on ${\cal T}_n$.\hfill
$\diamondsuit $

It is immediate to check that $\leq _B$ induces an order on the subset ${\cal Y}_n$ of ${\cal
T}_n$. This order is the order defined in [LRo].

\begin{rem} \label{wedg} Let $t=\bigvee (t^0,\dots ,t^k)$ and $z=\bigvee (z^0,\dots ,z^k)$ in ${\cal T}_n$, such that $\mid
t^i\mid = \mid z^i\mid $, for $0\leq i\leq k$.
\begin{enumerate}\item If $t\subseteq z$ (respectively, $t\leq _B z$), then $t^i\subseteq z^i$ (respectively
$t^i\leq _B z^i$), for all $0\leq i\leq k$.
\item If $w\in {\cal T}_n$ is such that $t\subseteq w\subseteq z$ (respectively $t\leq _B w \leq _B z$), then

\noi $w=\bigvee (w^0,\dots ,w^k)$ with $t^i\subseteq w^i\subseteq z^i$ (respectively $t^i\leq _Bw^i\leq _B z^i$),
for $0\leq i\leq k$.
\end{enumerate}
\end{rem}
\bs

\subsection{ Poset morphism from ${\cal P}_{\infty }$ to ${\cal T}_{\infty }$}

\begin{defn} Let us define a map $\Gamma : {\cal P}_{\infty } \longrightarrow {\cal T}_{\infty }$ as follows:
\begin{enumerate}
\item
$\Gamma ((0))=\downarrow $, where $\downarrow $ is the tree with one leaf and without vertex.
\item
Let $\gamma =\bigvee _{\omega }(\gamma _0,\dots ,\gamma _k)$ be an element of ${\cal P}_n$. The element
$\Gamma (\gamma )$ is the tree:
$$\Gamma (\gamma ):=\bigvee (\Gamma (\gamma _0),\dots ,\Gamma (\gamma _k)).$$
\end{enumerate}
\end{defn}

\begin{rem}\begin{enumerate}\item  It is immediate to see that $\Gamma ((1,\dots ,1))=\varsigma _n$, for $n\geq 1$.
\item An easy recursive argument shows that $\Gamma $ is surjective.
\item Let $\gamma _0,\dots ,\gamma _k$ be a family of elements in ${\cal P}_{\infty }$, with $\gamma _i\in {\cal
P}_{n_i,r_i}$ and $r=\sum _{i=0}^kr_i$, and let $\omega $ be an element in $SH(r_0,\dots ,r_k)$.
The images of the elements

\noi $\bigvee _{\omega }(\gamma _0,\dots ,\gamma _k)$ and $\bigvee _{1_r}(\gamma _0,\dots ,\gamma _k)$ under $\Gamma $ are
equal.
\item The image of the subset ${\cal P}_{n,r}$ of ${\cal P}_n$ under $\Gamma $ is the subset ${\cal T}_{n,r}:=\{ t\in {\cal
T}_n\quad :\quad t\quad {\rm has}\quad n-r\quad {\rm vertices}\} $.
\end{enumerate}
\end{rem}

\begin{thm}\label{presor} The map $\Gamma :{\cal P}_{\infty }\rightarrow {\cal T}_{\infty }$ verifies that:
\begin{enumerate} \item If $\gamma \subseteq \delta $ in ${\cal P}_n$, then $\Gamma (\gamma )\subseteq \Gamma (\delta )$
 in ${\cal T}_n$.
\item If $\gamma \leq _B \delta $ in ${\cal P}_n$, then $\Gamma (\gamma )\leq _B \Gamma (\delta )$
 in ${\cal T}_n$.
\end{enumerate}\end{thm}

\noi {\it Proof.} In this proof, given a sequence of integers $r_0,\dots ,r_k$, we shall denote by
$R_l$ the sum $R_l:=r_0+\dots +r_l$, for $0\leq l\leq k$.
\begin{enumerate}\item We know that $\gamma \subseteq \delta $ if, and only if, there exists a non-decreasing map
$\rho $ such that $\delta =\rho \circ \gamma $.
Since $\rho $ is a product of $t_i$\rq s, it suffices to prove the result for $\rho =t_j$.

The result is immediate for $n=0,1,2$. For $n\geq 3$, suppose $\gamma =\bigvee _{\omega }(\gamma _0,\dots ,\gamma _k)$, with
$\gamma _i\in {\cal P}_{n_i,r_i}$. For $j<\sum _{i=0}^kr_i$, we have that:
\begin{enumerate}\item If $t_j\circ \omega \in SH(r_0,\dots ,r_k)$, then $\delta =\bigvee _{t_j\circ \omega }(\gamma _0,\dots
,\gamma _k)$ and $\Gamma (\delta )=\Gamma (\gamma )$.
\item If $t_j\circ \omega \notin SH(r_0,\dots ,r_k)$, then Lemma \ref{prodwedge} implies that $\delta =\bigvee _{\omega
_j}(\gamma _0^j,\dots ,\gamma _k^j)$, where $\gamma _i^j=\gamma _i$ or $\gamma _i^j=t_{j_i}\circ \gamma _i$. So, $\gamma
_i\subseteq \gamma _i^j$, for $0\leq i\leq k$.

\noi By recursive hypothesis $\Gamma (\gamma _i)\subseteq \Gamma (\gamma _i^j)$, and
$\Gamma (\delta )=\bigvee (\Gamma (\gamma
_0^j),\dots ,\Gamma (\gamma _k^j))$. So, we get that
$\Gamma (\gamma )\subseteq \Gamma (\delta )$.
\end{enumerate}
For $j=\sum _{i=0}^kr_i=r-1$, the element $\delta $ is obtained as a wedge of two types of elements:
$$\left \{ \begin{array}{cc}
\gamma _i&\quad {\rm for}\quad \omega (R_i)< r-1,\\
\gamma _{l,i}&\quad {\rm for}\quad \omega (R_i)= r-1\quad {\rm and}\quad \gamma _i=\bigvee _{\omega ^i}(\gamma
_{0,i},\dots ,\gamma _{k_i,i}).\end{array}\right .$$
So, $\Gamma (\delta )$ is the wedge of the elements:
$$\left \{ \begin{array}{cc}
\Gamma (\gamma _i)&\quad {\rm for}\quad \omega (R_i)< r-1,\\
\Gamma (\gamma _{l,i})&\quad {\rm for}\quad \omega (R_i)= r-1\quad {\rm and}\quad \gamma _i=\bigvee _{\omega
^i}(\gamma _{0,i},\dots ,\gamma _{k_i,i}).\end{array}\right .$$
The argument above implies that $\Gamma (\delta )$ is obtained from $\Gamma (\gamma )$ by contracting all the edges which
joint the root of $\Gamma (\gamma _i)$ to the root of $\Gamma (\gamma )$, for all $0\leq i\leq k$ such that $\omega
(R_i)= r-1$. It implies that $\Gamma (\gamma )\subseteq \Gamma (\delta )$.

\item For the weak Bruhat order the result is immediate for $n\leq 2$. For $n\geq 3$, suppose
$\gamma =\bigvee _{\omega }(\gamma _0,\dots ,\gamma _k)$, and $\delta =\bigvee _{\upsilon }(\delta
_0,\dots ,\delta _h)$, with $\gamma _i\in {\cal P}_{n_i,r_i}$ and $\delta _i\in {\cal P}_{m_i,s_i}$.

\noi a) If $\delta =t_j\circ \gamma $, with $\gamma ^{-1}(j)<\gamma ^{-1}(j+1)$, then we have that:

\noi $(i)$ For $j<r-1$ and $t_j\circ \omega \in SH(r_0,\dots ,r_k)$, it holds that
$\Gamma (\delta )=\Gamma (\gamma )$.

\noi $(ii)$ For $j<r-1$ and $t_j\circ \omega \notin SH(r_0,\dots ,r_k)$, Lemma \ref{prodwedge} states that

\noi $\delta =\bigvee _{\omega _j}(\gamma _0^j,\dots ,\gamma _k^j)$, with
$$\left \{ \begin{array}{cc}
\gamma _i^j=\gamma _i&\quad {\rm when}\quad \{ j,j+1\} \not\subseteq \omega
(\{ R_{i-1}+1, \dots ,R_i\} ),\\
\gamma _i^j=t_{j_i}\circ \gamma _i&\quad {\rm when}\quad \{ j,j+1\} \subseteq \omega (\{ R_{i-1}+1, \dots
,R_i\} ).\end{array}\right .$$
The fact that $\gamma ^{-1}(j)<\gamma ^{-1}(j+1)$ implies that there exists at most one $i_0$ such that $\gamma
_{i_0}^j\not = \gamma _{i_0}$ and $\gamma _{i_0}^{-1}(j_{i_0})<\gamma _{i_0}^{-1}(j_{i_0}+1)$.
So,
$$\Gamma (\delta )=\bigvee (\Gamma (\gamma _0),\dots ,\Gamma (t_{j_{i_0}}\circ \gamma _{i_0}),\dots ,
\Gamma (\gamma _k)),$$
 and, by recursive hypothesis, $\Gamma (\gamma _{i_0})\leq _B \Gamma
(t_{j_{i_0}}\circ \gamma _{i_0})$, which implies that

\noi $\Gamma (\gamma )\leq _B\Gamma (\delta )$.

\noi $(iii)$ For $j=r-1$, we get that $\gamma ^{-1}(r-1)\subseteq \{ 1,\dots ,r_0\} $.
It holds that
$$\Gamma (\gamma )=\bigvee (\Gamma (\gamma _0),\dots ,\Gamma (\gamma _k))$$ and that $\Gamma (\delta )$ is the wedge of
$\Gamma (\gamma _{0,0}),\dots ,\Gamma (\gamma _{k_0,0}),\Gamma (\gamma _1),\dots ,\Gamma (\gamma _k)$, with

\noi $\Gamma (\gamma _0)=\bigvee (\Gamma (\gamma _{0,0}),\dots ,\Gamma (\gamma _{k_0,0}))$.
So, we get that $\Gamma (\gamma )\leq _B \Gamma (\delta )$.

\noi b) The proof of the case $t_j\circ \delta =\gamma $, with $\delta ^{-1}(j)>\delta ^{-1}(j+1)$, can be obtained following
the same arguments that we use to prove a). \hfill $\diamondsuit $
\end{enumerate}

In order to see that the orders $\subseteq $ and $\leq _B$ of ${\cal P}_{\infty }$ induce the orders $\subseteq $ and $\leq
_B$ on ${\cal T}_{\infty }$, we need to prove some additional results.

\begin{rem} Let $\gamma $, $\omega $ and $\delta $ be elements of ${\cal P}_m$ such that $\Gamma (\gamma )=\Gamma (\delta )$
and $\gamma \subseteq \omega \subseteq \delta $. The theorem above implies that $\Gamma (\omega )=\Gamma (\gamma )=\Gamma
(\delta )$.
\end{rem}

Let $\leq $ be the relation on ${\cal T}_{\infty }$ given by $t\leq z$ if and only if there exist $\gamma \in \Gamma ^{-1}(t)$
and $\delta \in \Gamma ^{-1}(z)$ such that $\gamma
\subseteq \delta $.

\begin{prop} The relation $\leq $ defines a partial order on ${\cal T}_{\infty }$ which coincides with $\subseteq $.
\end{prop}

\noi {\it Proof.} It results easy to check that $\leq $ is reflexive and transitive. Lemma \ref{prodwedge} states that if
$\delta =\rho \circ \gamma $, with $\rho $ a non decreasing element in ${\cal P}_r$, then either $\Gamma (\delta )=\Gamma
(\gamma )$ or $\Gamma (\delta )$ is obtained from $\Gamma (\gamma )$ by contracting some internal edges. This fact implies the
antisymmetry of $\leq $.

If $t\leq z$, then there exist $\gamma \in \Gamma ^{-1}(t)$ and $\delta \in \Gamma ^{-1}(z)$ with
$\gamma \subseteq \delta $. By Theorem \ref{presor}, it holds that $t=\Gamma (\gamma )\subseteq \Gamma (\delta )=z$.

If $t\subseteq z$, then is becomes clear that $t\leq z$ for $\mid t\mid =\mid z\mid \leq 2 $. For $\mid t\mid =\mid z\mid >2$,
suppose that $z$ is obtained from $t=\bigvee (t^0,\dots ,t^k)$ by contracting an internal edge.
\begin{enumerate}\item If $z=\bigvee (t^0,\dots ,z^j,\dots ,t^k)$, with $t^j\subseteq z^j$, then by recursive hypothesis
there exist $\gamma _j\in \Gamma ^{-1}(t^j)$ and $\delta _j\in \Gamma ^{-1}(z^j)$ such that $\delta _j=\rho _j\circ
\gamma _j$, with $\rho _j$ a non-decreasing map in ${\cal P}_{\infty }$.

\noi Choose an element $\gamma _i\in \Gamma ^{-1}(t^i)$, for $i\not = j$. The element $\gamma :=\bigvee_{1_{r-1}}(\gamma _0,
\dots ,\gamma _k)$ belongs to $\Gamma ^{-1}(t)$. Consider $\delta := \bigvee _{1_{r-1-r_j+s_j}}(\gamma _0,\dots ,\delta
_j,\dots ,\gamma _k)$, we get that $\delta \in \Gamma ^{-1}(z)$ and $\gamma \subseteq \delta $, which implies $t\leq z$.

\item For $z=\bigvee (t^0,\dots ,t^{0,j},\dots ,t^{n_j,j},\dots ,t^k)$, with $t^j=\bigvee (t^{0,j},\dots ,t^{n_j,j})$, choose
elements $\gamma _i$ be in $\Gamma ^{-1}(t^i)$, for $i\not = j$, and $\gamma _{l,j}\in \Gamma ^{-1}(t^{l,j})$, for $0\leq l\leq
n_j$.

\noi The element $\gamma _j:=\bigvee _{1_{r_j-1}}(\gamma _{0,j},\dots ,\gamma _{n_j,j})$ belongs to $\Gamma ^{-1}(\gamma _j)$,
$\gamma :=\bigvee _{1_{r-1}}(\gamma _0,\dots ,\gamma _k)$ is in $\Gamma ^{-1}(t)$, and $\delta :=\bigvee
_{1_{r-1-r_j+\sum r_{l,j}}}(\gamma _0,\dots ,\gamma _{0,j},\dots ,\gamma _{n_j,j},\dots ,\gamma _k)$ belongs to $\Gamma
^{-1}(z)$. Since $\gamma \subseteq \delta $ we may assert that $t\leq z$.
\hfill $\diamondsuit $
\end{enumerate}

\begin{defn} Given a tree $t\in {\cal T}_m$, we define the elements $Min(t)$ and $Max(t)$ in ${\cal P}_m$ recursively, as
follows:
\begin{enumerate}\item $Min(\downarrow ):=(0)=:Max(\downarrow )$.
\item For $t=\bigvee (t^0,\dots ,t^k)$, with $t^i\in {\cal T}_{n_i}$ and $n=\sum _{i=0}^kn_i$, the degree of $t$ is $n+k$.
Define
$$\displaylines {
Min(t):=\bigvee _{1_{r-1}}(Min(t^0),\dots ,Min(t^k)),\cr
Max(t):=\bigvee _{\xi _{s_0,\dots ,s_k}}(Max(t^0),\dots ,Max(t^k)),\cr }$$
where $\xi _{s_0,\dots ,s_k}$ is the biggest element of $SH(s_0,\dots ,s_k)$ for the weak Bruhat order, $Min(t^i)\in {\cal
P}_{n_i,r_i}$, $Max(t^i)\in {\cal P}_{n_i,s_i}$ and $\sum _{i=0}^kr_i=r-1$.
\end{enumerate}
\end{defn}

\begin{prop} Let $t\in {\cal T}_m$ be a tree, the inverse image under $\Gamma $ of $t$ is an interval for the weak Bruhat
order. Moreover, it verifies
$$\Gamma ^{-1}(t)=\{ \delta \in {\cal P}_m\quad {\rm such\quad that}\quad Min(t)\leq _B \delta \leq _B Max(t)\} .$$
\end{prop}

\noi{\it Proof.} The assertion is obviously true for $m=0,1,2$. For $m\geq 3$, let $t\in {\cal T}_m$ with $t=\bigvee
(t^0,\dots ,t^k)$. By recursive hypothesis, we know that
$$\Gamma ^{-1}(t^i)=\{  \delta _i\in {\cal P}_{n_i}\quad {\rm such\quad that}\quad Min(t^i)\leq _B \delta _i\leq _B Max(t^i)\}
.$$

Suppose $\gamma =\bigvee _{\omega }(\gamma _0,\dots ,\gamma _k)$ is such that $\Gamma (\gamma )=t$.  We have that $\Gamma
(\gamma _i)=t^i$, for $0\leq i\leq k$. By Corollary \ref{wBor1} and Lemma \ref{wBor2} we obtain that:
$$Min(t)=\bigvee _{1_{r-1}}(Min(t^0),\dots ,Min(t^k))\leq _B \bigvee _{1_{r\rq -1}}(\gamma _0,\dots ,\gamma _k)\leq _B \gamma
,$$
and, there exists $\omega \rq \leq _B \xi _{s_0,\dots ,s_k}$ such that:
$$\gamma \leq _B\bigvee _{\omega \rq }(Max(t^0),\dots ,Max(t^k))\leq _B
\bigvee _{\xi _{s_0,\dots ,s_k}}(Max(t^0),\dots ,Max(t^k))=Max(t).$$

Conversely, if $Min(t)\leq _B \gamma \leq _B Max(t)$, then Theorem \ref{presor} states that
$$t=\Gamma (Min(t))\leq _B \Gamma (\gamma )\leq _B \Gamma (Max(t))=t.\hfill \diamondsuit $$

So, the weak Bruhat order of ${\cal P}_{\infty }$ induces a partial order on ${\cal T}_{\infty }$.

\begin{prop} The weak Bruhat order $\leq _B$ of ${\cal P}_{\infty }$ induces the weak Bruhat order on ${\cal T}_{\infty }$.
\end{prop}

\noi {\it Proof.} Using that $\Gamma $ preserves the weak Bruhat order, it is easily seen that $t\leq z$ for the order
induced by the weak Bruhat order of ${\cal P}_{\infty }$ implies that $t\leq _B z$.

It is immediate to check that $t\leq _B z$ implies $t\leq z$ for the induced order, when $\mid t\mid =\mid z\mid \leq 2$.
For $\mid t\mid \geq 3$, we have to consider the following three cases:
\begin{enumerate}\item If $t=\bigvee (t^0,\dots ,t^k)$ and $z=\bigvee (z^0,\dots ,z^k)$, with $t^i\leq _B z^i$ for all $i$,
then by recursive hypothesis we get that $t^i\leq z^i$ for the induced order.

\noi So, there exist two families $\gamma _0,\dots
,\gamma _k$ and $\delta _0,\dots ,\delta _k$ of element in ${\cal P}_{\infty }$ such that $\gamma _i\leq _B\delta _i$ ,
$\Gamma (\gamma _i)=t^i$ and $\Gamma (\delta _i)=z^i$.
Since $t=\Gamma (\bigvee _{1_{r-1}}(\gamma _0,\dots ,\gamma _k))$,
$z=\Gamma (\bigvee _{1_{s-1}}(\delta _0,\dots ,\delta _k))$ and $\bigvee _{1_{r-1}}(\gamma _0,\dots ,\gamma _k)\leq _B \bigvee
_{1_{s-1}}(\delta _0,\dots ,\delta _k)$, we get that $t\leq z $ for the induced order.

\item If $t=\bigvee (t^0,\dots ,t^k)$ and $z=\bigvee (t^{0,0},\dots ,t^{k_0,0},t^1,\dots ,t^k)$, for $t^0=\bigvee
(t^{0,0},\dots ,t^{k_0,0})$, then there exist $\gamma _{0,0},\dots ,\gamma _{k_0,0},\gamma _1,\dots ,\gamma _k$ in ${\cal
P}_{\infty }$ such that $\Gamma (\gamma _{l,0})=t^{l,0}$, for $0\leq l\leq k_0$, and $\Gamma (\gamma _i)=t^i$, for $1\leq i\leq
k$.

\noi Consider the elements $\gamma _0 :=\bigvee _{1_{n_0}}(\gamma _{0,0},\dots ,\gamma _{k_0,0})$, $\gamma :=\bigvee
_{\xi _{r_0,r-1-r_0}}(\gamma _0,\dots ,\gamma _k)$ and $\delta :=\bigvee _{1_{s-1}}(\gamma _{0,0},\dots ,\gamma _{k_0,0},
\gamma _1,\dots ,\gamma _k)$,  where $\gamma _i\in {\cal P}_{n_i,r_i}$, $\gamma _{l,0}\in {\cal P}_{n_{l,0},r_{l,0}}$, $\sum
_{i=0}^kr_i=r-1$ and $\sum _{l=0}^{k_0}r_{l,0}+\sum _{i=1}^kr_i=s-1$. It holds that $\Gamma (\gamma )=t$, $\Gamma (\delta
)=z$, $\gamma ^{-1}(r-1)<\gamma ^{-1}(r)$ and $\delta = t_{r-1}\circ \gamma $, which implies $t\leq z $ for the induced order.

\item The proof for the case $t=\bigvee (z^0,\dots ,z^k,z^{0,k},\dots ,z^{0,h_k})$ and $z=\bigvee (z^0,\dots ,z^k)$ is obtained
in the same way that the one of case $2$. \hfill $\diamondsuit $
\end{enumerate}
\bs

\subsection{Free dendriform trialgebra structure on $k[{\cal T}_{\infty }]$}

Note that the operations $\succ $, $\cdot $ and $\prec $ defined in Section 2 induce operations on the $k$-vector space
$k[{\cal T}_{\infty }]$, spanned by the set of planar trees. The operations are defined as follows:
\begin{enumerate}
\item $t\succ z=\sum w$, where the sum is taken over the set
$\{ w\in {\cal T}_{\infty }$ such that there exist $\gamma \in \Gamma ^{-1}(t),\quad \delta \in \Gamma ^{-1}(z)$ and $\sigma
\in SH^{\succ }(r,s)$ such that $\sigma \circ (\gamma \times \delta )\in \Gamma ^{-1}(w)\} .$
\item $t\cdot z=\sum w$, where the sum is taken over the set
$\{ w\in {\cal T}_{\infty }$ such that there exist $\gamma \in \Gamma ^{-1}(t),\quad \delta \in \Gamma ^{-1}(z)$ and
$\sigma \in SH^{\bullet }(r,s)$ such that $\sigma \circ (\gamma \times \delta )\in \Gamma ^{-1}(w)\} .$
\item $t\prec z=\sum w$, where the sum is taken over the set $\{ w\in {\cal T}_{\infty }$ such that there
exist $\gamma \in \Gamma ^{-1}(t),\quad \delta \in \Gamma ^{-1}(z)$ and $\sigma
\in SH^{\prec }(r,s)$ such that $\sigma \circ (\gamma \times \delta )\in \Gamma ^{-1}(w)\} .$
\end{enumerate}

\begin{rem}\label{varios} Given a partition $n_0,\dots ,n_k$ of $n$ and $1\leq j\leq n$, we denote by $SH(n_0,\dots ,n_k)_r$
the set of elements $\omega \in SH(n_0,\dots ,n_k)\cap {\cal P}_{n,r}$. An easy calculation shows that the following
equalities hold:
$$\displaylines {
SH(r,s_0,\dots ,s_h)=\bigcup _{j} SH(r,j)\circ (1_r\times SH(s_0,\dots ,s_h))=\hfill \cr
\hfill \bigcup _{j}SH(j,s_1,\dots ,s_h)\circ (SH(r,s_0)\times 1_{s-s_0-1}), \qquad {\bf (1)}\cr
\bigcup _{i,j}SH(i,j)\circ (SH(r_0,\dots ,r_k)_i\times SH(s_0,\dots ,s_h)_j)=\hfill \cr
\hfill \bigcup _{i}SH(r_0,\dots ,r_{k-1},i,s_1,\dots ,s_h)\circ (1_{r-r_k-1}\times SH(r_k,s_0)\times
1_{s-s_0-1}), \qquad {\bf (2)}\cr
SH(r_0,\dots ,r_k,s)=\bigcup _{j}SH(j,s)\circ (SH(r_0,\dots ,r_k)\times 1_s)=\hfill \cr
\hfill \bigcup _{j}SH(r_0,\dots ,r_{k-1},j)\circ (1_{r-r_k-1}\times SH(r_k,s)), \qquad {\bf (3)}\cr }$$
where $r-1=\sum _{i=0}^kr_i$ and $s-1=\sum _{i=0}^hs_i$.
\end{rem}

\begin{thm}\label{qfin} The operations $\succ $, $\cdot $ and $\prec $ induced on $k[{\cal T}_{\infty }]$ by the dendriform
trialgebra structure of $k[{\cal P}_{\infty }]$ verify the following equalities:
\begin{enumerate}\item $t*\downarrow=t=\downarrow *t$, for all $t\in {\cal T}_{\infty }$,
\item $t\succ z=\bigvee (t*z^0,z^1,\dots ,z^h)$,
\item $t\cdot z=\bigvee (t^0,\dots ,t^{k-1},t^k*z^0,z^1,\dots ,z^h)$,
\item $t\prec z = \bigvee (t^0,\dots ,t^{k-1},t^k*z)$,
\end{enumerate}
for $t=\bigvee (t^0,\dots ,t^k)$ and $z=\bigvee (z^0,\dots ,z^h)$.
\end{thm}

\noi {\it Proof.} $1$ is obvious.

The proofs of $2$, $3$ and $4$ use similar arguments, applying the Remark above. We give the detailed proof of $3$ (which
seems to be just a bit more difficult), the other ones may be obtained following analogous steps.

$3$. We have that $SH^{\bullet }(j,l)=(SH(j-1,l-1)\times 1_1)\circ z(j-1,l-1)$. A tree $w$ appears with coefficient $1$ in
 $t\cdot z$ if there  exist $\{ \gamma _i\in {\cal P}_{n_i,r_i}\} _{0\leq i\leq k}$, $\omega \in SH(r_0,\dots ,r_k)_{j-1}$,
$\{ \delta _i\in {\cal P}_{m_i,s_i}\}_{0\leq i\leq h}$, $\upsilon \in SH(s_0,\dots ,s_h)_{l-1}$ and $\sigma \in SH(j-1,l-1)$
such that the element
$\sigma \circ (\gamma \times \delta)$ belongs to $\Gamma ^{-1}(w)$, with $\sigma =(\sigma _1\times 1_1)$,
 $\gamma =\bigvee _{\omega }(\gamma _0\times \dots \gamma _k)$ and $\delta =\bigvee _{\upsilon }(\delta _0\times \dots
\times \delta _h)$. Using Remark \ref{varios} {\bf (2)}, we get that
$$\displaylines {
\sigma \circ (\gamma \times \delta)= (\sigma _1\times 1_1)\circ z(j-1,l-1)\circ \hfill \cr
\hfill ((\omega \circ (\gamma _0\times \dots
\gamma _k))\times 1_1\times  (\upsilon \circ (\delta _0\times \dots \times \delta _h))\times 1_1)\circ \hfill \cr
\hfill (z(n_0,\dots ,n_k)\times z(m_0,\dots ,m_h))=\cr
(\sigma _1\times 1_1)\circ ((\omega \circ (\gamma _0\times \dots
\gamma _k)\times \upsilon \circ (\delta _0\times \dots \times \delta _h)\times 1_1)\hfill \cr
\hfill \circ z(n_0,\dots ,n_k,m_0,\dots ,m_h)=\cr
 \tau \circ (1_{r-r_0-1}\times \zeta \times 1_{s-s_0-1})\circ \hfill \cr
\hfill (\gamma _0\times \dots \times \gamma _k\times \delta _0\times \dots \times \delta _h\times 1_1)\circ z(n_0,\dots
,n_k,m_0,\dots ,m_h)=\cr
\hfill \bigvee _{\tau }(\gamma _0,\dots ,\gamma _{k-1},\zeta \circ (\gamma _k\times \delta _0),\delta _1,\dots ,\delta _h),\cr
 }$$
for unique elements $\zeta \in SH(r_k,s_0)_i$ and $\tau \in SH(r_0,\dots ,r_{k-1},i,s_1,\dots ,s_h)$. So, $w$ appears in
$\bigvee (t^0,\dots ,t^{k-1},t^k*z^0,z^1,\dots ,z^h)$.

Looking at the same equalities from the lower line to the upper one , we show that any tree $w$ which appears
in $\bigvee (t^0,\dots ,t^{k-1},t^k*z^0,z^1,\dots ,z^h)$ with coefficient $1$ appears with the same coefficient in $t\cdot
z$, too. \hfill $\diamondsuit $

The following result is a consequence of the Theorem above and the description of the free dendriform trialgebra given in
\cite{LR4}.

\begin{cor} The dendriform trialgebra $(k[{\cal T}_{\infty }],\succ ,\cdot ,\prec )$ is the free dendriform trialgebra
spanned by one element.
\end{cor}

Let us describe now the operations of the free dendriform algebra $k[{\cal T}_{\infty }]$ in terms of the weak Bruhat order of
${\cal T}_{\infty }$.

Given a tree $t\in {\cal T}_n$, we numerate its leaves from left to right by $0,1,\dots n$.

\begin{notn} Let $t\in {\cal T}_n$ and $z\in {\cal T}_m$ be two planar trees.
\begin{enumerate} \item For any $0\leq j\leq m$, the grafting
$t\circ _j z$ of t on the $j$-th leaf of $z$ is the tree obtained by joining the root of $t$ to the leaf $j$ of $z$.
\item We denote by $t/z $ the tree $t\circ _0 z$, and by $t\backslash z$ the tree $z\circ _n t$.
\end{enumerate}
\end{notn}

\begin{prop}\label{fin} Let $t\in {\cal T}_n$ and $z\in {\cal T}_m$. Given elements $\gamma \in \Gamma ^{-1}(t)\cap {\cal
P}_{n,r}$ and
$\delta \in \Gamma ^{-1}(z)\cap {\cal P}_{m,s}$, it holds that:
\begin{enumerate} \item $\Gamma (\gamma \times \delta )=t/ z,$
\item $\Gamma (\xi _{r,s}\circ (\gamma \times \delta ))=t\backslash z$.
\end{enumerate}
\end{prop}

\noi {\it Proof.} \begin{enumerate}\item If $z =\downarrow $, then $\gamma \times (0) = \gamma $ and
$\Gamma (\gamma )=t=t/ \downarrow $.

\noi If $m\geq 1$, then $\delta =\bigvee _{\upsilon }(\delta _0,\dots ,\delta _h)$, for a family $\delta _i\in \Gamma
^{-1}(z^i)$, with $0\leq i\leq h$, and $\upsilon \in SH(s_0,\dots ,s_h)$. We have that:
$$\displaylines {
\gamma \times \delta =(1_r\times \upsilon )\circ (\gamma \times \delta _0\times \dots \times \delta _h\times
1_1) \circ (1_n\times z(m_0,\dots ,m_h)=\cr
(1_r\times \upsilon )\circ ((\gamma \times \delta _0)\times \delta _1\times \circ \times \delta _h\times 1_1)\circ
z(n+m_0,m_1,\dots ,m_h)=\cr
\bigvee _{1_r\times \upsilon}(\gamma \times \delta _0,\delta _1,\dots ,\delta _h).\cr }$$
The recursive hypothesis states that $\Gamma (\gamma \times \delta _0)=t/z^0$, which implies that
$\Gamma (\gamma \times \delta )=\bigvee (t/z^0,z^1,\dots ,z^h)=t/z$.

\item If $t=\downarrow $, then $\xi _{0,r}\circ ((0)\times \delta )=\delta $ and $\Gamma (\delta )=z=\downarrow \backslash z$.

\noi If $n\geq 1$, then $\gamma =\bigvee _{\omega }(\gamma _0,\dots ,\gamma _k)$, with $\gamma _i\in \Gamma ^{-1}(t^i)$ and
$\omega \in SH(r_0,\dots ,r_k)_r$, with $\sum _{i=0}^kr_i=r\rq -1$. It holds that
$$\displaylines {
\xi _{r,s}\circ (\gamma \times \delta )=\xi _{r,s}\circ (\omega \times 1_{s+1})\circ (\gamma _0\times \dots \times \gamma
_k\times 1_1\times \delta )\circ (z(r_0,\dots ,r_k)\times 1_m)=\cr
(\omega (1)+s,\dots ,\omega (r\rq -1)+s,r+s,1,\dots ,s)\circ (\gamma _0\times \dots \times \gamma
_k\times 1_1\times \delta )\circ (z(r_0,\dots ,r_k)\times 1_m)=\cr
((\xi _{r-1,s}\circ (\omega \times 1_s)\circ (\gamma _0\times \dots \times \gamma _{k-1}\times (\gamma
_k\times \delta ))\times 1_1)\circ z(r_0,\dots ,r_{k-1},r_k+s)=\cr
(\upsilon \times 1_1)\circ
(\gamma _0\times \dots \times \gamma _{k-1}\times (\xi _{r_k,s}\circ (\gamma _k\times \delta ))\times 1_1)\circ z(r_0,\dots
,r_{k-1},r_k+s)=\cr
\bigvee _{\upsilon}(\gamma _0,\dots ,\gamma _{k-1},(\xi _{r_k,s}\circ (\gamma _k\times \delta ))),\cr }$$
where $\upsilon :=$
$$(\omega (1)+s,\dots ,\omega (r_0+\dots +r_{k-1})+s,1,\dots ,s,\omega (r_0+\dots +r_{k-1}+1)+s,\dots ,\omega
(r\rq -1)+s),$$
belongs to $SH(r_0,\dots ,r_{k-1},r_k+s)$.

\noi By recursive hypothesis, $\Gamma (\xi _{r_k,s}\circ (\gamma _k\times \delta )=t^k\backslash z$.

\noi So, we get that $\Gamma (\xi _{r,s}\circ (\gamma \times \delta ))=\bigvee (t^0,\dots ,t^{k-1},t^k\backslash z)
=t\backslash z$. \hfill $\diamondsuit $
\end{enumerate}

Theorem \ref{qfin} and Proposition \ref{fin} imply the following result:

\begin{cor} The operations $\succ $, $\cdot $, $\prec $ and $*$ of the free dendriform trialgebra $k[{\cal T}_{\infty }]$ are
defined in terms of the weak Bruhat order and the products $/$ and $\backslash $ by the following conditions:
\begin{enumerate}\item $t\succ z =\sum _{t/z\leq _B w\leq _B \bigvee (t\backslash z^0,z^1,\dots ,z^k)}w,$
\item $t\cdot z=\sum _{\bigvee (t^0,\dots ,t^{k-1},t^k/z^0,z^1,\dots ,z^h)\leq _B w\leq _B \bigvee  (t^0,\dots
,t^{k-1},t^k\backslash z^0,z^1,\dots ,z^h)}w,$
\item $t\prec w =\sum _{\bigvee (t^0,\dots ,t^{k-1},t^k/z)\leq _B w \leq _B t\backslash z}w,$
\item $t*z=\sum _{t/z\leq _B w\leq _B t\backslash z}w.$
\end{enumerate}
\end{cor}
\bs

\noi {\bf Final comment} Given a Coxeter system $(W,S)$ there exists another order on the set ${\cal P}^{(W,S)}$ which extends
the weak Bruhat order of $W$. Consider the transitive relation $\leq _C$ spanned by:
$$ W_J\circ w \leq _C W_J\circ (\alpha _{J,s}^{-1}\circ w),\ {\rm for}\ J\subseteq S,\ s\in S\setminus J\ {\rm and}\ w\in
X_{J\cup \{ s\} }^{-1},$$
where $W_J$, $X_K$ and $\alpha _{J,s}$ are the elements defined in Section 1.

\noi Observe that if $x\leq _C y$, then $x\leq _B y$. So, the weak Bruhat order is weaker than $\leq _C$.

\noi This new order preserves the graduation of ${\cal P}^{(W,S)}$, that is two elements $x$ and $y$ of ${\cal P}^{(W,S)}$
are comparable for $\leq _C$ only if $x,y\in {\cal P}_r^{(W,S)}$, for some $0\leq r\leq \mid S\mid $.

\noi Mimicking the results of Section 3 it is possible to see that the order $\leq _C$ of ${\cal P}_{\infty }$ induces, via
the map $\Gamma $, a partial order on ${\cal T}_{\infty }$.

\noi The order $\leq _C$ on ${\cal P}_n$ is the transitive relation generated by:
$$omega\leq _Cs_i\circ \omega,\quadÊ{\rm for}\quad \omega ^{-1}(i)\leq \omega ^{-1}(i+1),$$
where $\omega ^{-1}(j)$ denotes the set of elements $l$ of $\{ 1,\dots ,n\}$ such that $\omega (l)=j$.

\noi On ${\cal T}_n$, the order $\leq _C$ is transitively spanned by the relations:
\begin{enumerate} \item If $t^{i_0} \leq _C w^{i_0} \in {\cal T}_{n_{i_0}}$ , then
$\bigvee(t^0, \ldots , t^{i_0},\ldots , t^k) \leq _C \bigvee(t^0, \ldots , w^{i_0},\ldots , t^k)$ in
${\cal T}_{n_0 + \ldots + n_k + k -1}$,
\item If $t=\bigvee(t^0, \ldots, t^k)\in {\cal T}_n$ and $w=\bigvee(w^0,
\ldots, w^h) \in {\cal T}_m$, then
$$\bigvee(t,w^0, \ldots, w^h) \leq _C \bigvee(t^0, \ldots, t^k, w).$$
\end{enumerate}

In this case, the order $\leq _C$ also permits to construct associative algebra structures on $k[{\cal P}_{\infty }]$ and
$k[{\cal T}_{\infty }]$. These structures coincide with the associative algebra structures defined on $k[{\cal P}_{\infty }]$
and $k[{\cal T}_{\infty }]$ by F. Chapoton in \cite{Ch}.

\end{document}